\documentclass[a4paper,10pt,english]{article}

\usepackage[latin1]{inputenc}

\usepackage{calc}
\usepackage{graphicx,palatino,verbatim}

\usepackage{amsmath}

\usepackage[T1]{fontenc}
\usepackage[english]{babel}
\usepackage{newlfont}

\usepackage{amsthm}

\usepackage{amssymb}

\usepackage{amsmath}

\usepackage{eufrak}

\usepackage{latexsym}

\usepackage{young}
\usepackage[mathscr]{eucal}

\usepackage{amsfonts}
\usepackage{syntonly}
\usepackage{graphicx}
\usepackage{mathptmx}
\usepackage{subfig}
\usepackage{url}
\usepackage[letterpaper]{geometry}
\input xy
\xyoption{all}

\newtheorem{teo}{Theorem}

\newtheorem{teorema}{Theorem}[section]

\newtheorem{oss}[teorema]{Remark}

\newtheorem{lemma}[teorema]{Lemma}
\newtheorem{definizione}[teorema]{Definition}
\newtheorem{esempi}[teorema]{Example}

\newtheorem{proposizione}[teorema]{Proposition}

\newcommand{\keywords}[1]{%
  \begingroup 
  \renewcommand\subsubsection{}%
  \addtocounter{subsubsection}{-1}%
  \subsubsection{\small{\textbf{Keywords}:} #1}%
  \endgroup} 

\begin{document}
\title{Semi-invariants of symmetric quivers of tame type}
\author{Riccardo Aragona}
\date{}
\maketitle
\begin{center}
Università degli Studi di Roma \textquotedblleft
Tor Vergata\textquotedblright\\
Dipartimento di Matematica\\
Via della Ricerca Scientifica 1, 00133 Rome (Italy)
\end{center}

\begin{center}
E-mail: $\quad$ ric$\_$aragona@yahoo.it
\end{center}
\begin{abstract}
\maketitle A symmetric quiver $(Q,\sigma)$ is a 
 finite quiver without oriented cycles $Q=(Q_0,Q_1)$ equipped with a contravariant involution $\sigma$ on $Q_0\sqcup Q_1$. The
involution allows us to define a nondegenerate bilinear form $<,>$
on a representation $V$ of $Q$. We shall say that $V$ is orthogonal
if $<,>$ is symmetric and symplectic if $<,>$ is
skew-symmetric. Moreover, we define an action of products of
classical groups on the space of orthogonal representations and on
the space of symplectic representations. So we prove that if $(Q,\sigma)$ is a symmetric quiver of tame type then the rings of semi-invariants for this
action are spanned by the semi-invariants of determinantal type
$c^V$ and, when matrix defining $c^V$ is
skew-symmetric, by the Pfaffians $pf^V$. To prove it, moreover, we describe the symplectic and orthogonal generic decomposition of a symmetric dimension vector.
\end{abstract}
\begin{center}
\end{center}

\section*{Introduction}
 The symmetric quivers and their orthogonal and symplectic representations have been introduced by Derksen and Weyman in \cite{dw2} to provide a formalization in the quiver setting of some problem related to representations of classical groups. Generalizations of quivers and their representations have been defined in a different setting by Zubkov in \cite{z} and Shmelkin in \cite{sh}.\\
The representations of symmetric quivers could be a tool to classify products of flag varieties with finitely many orbits under the diagonal action of classical groups. Magyar, Weyman and Zelevinsky in \cite{mwz} solved this problem for general linear groups using usual quiver setting.\\
It would be also interesting to generalize to symmetric quivers the results about virtual representations and virtual semi-invariants of quivers given by Igusa, Orr, Todorov and Weyman in \cite{iotw}.\\ 
The author,  in \cite{ar}, displayed a set of generators of rings of semi-invariants of symmetric quivers of finite type. In this paper, the second one extract from his PhD thesis, supervised by Professor Jerzy Weyman,  we provide similar results for symmetric quivers of tame type. Similar problems, in a different setting, have been studied by Lopatin in \cite{l}, using ideas from \cite{lz}.\\
Analogous results for usual quivers have been obtained independently by Derksen and Weyman in \cite{dw1}, Domokos and Zubkov in \cite{dz} and Schofield and Van den Bergh in \cite{sv}. \\
Let $Q=(Q_0,Q_1)$ be  a quiver, where $Q_0$ and $Q_1$ are respectively the set of vertices and the set of arrows of $Q$, and let $\sigma$ be an involution on $Q_0\sqcup Q_1$. The pair $(Q,\sigma)$ is called \textit{symmetric quiver}.\\ 
Let $V$ be a representation of $Q$. The involution allows us to define
a nondegenerate bilinear form $<,>$ on $V$. We call the pair $(V,<,>)$ symplectic
(respectively orthogonal) representation of $(Q,\sigma)$ if $<,>$ is skew-symmetric
(respectively symmetric). We define $SpRep(Q,\beta)$ and
$ORep(Q,\beta)$ to be respectively the space of symplectic
$\beta$-dimensional representations and the space of orthogonal
$\beta$-dimensional representations  of $(Q,\sigma)$. Moreover we
can define an action of a product of classical groups, which we
call $SSp(Q,\beta)$ in the symplectic case and $SO(Q,\beta)$ in
the orthogonal case, on these space. Let $SpSI(Q,\beta)$ and $OSI(Q,\beta)$ be respectively the ring of symplectic semi-invariants and the ring of orthogonal semi-invariants of a symmetric quiver $(Q,\sigma)$.\\
Let $(Q,\sigma)$ be a symmetric quiver and $V$ a representation of
the underlying quiver $Q$ such that
  $\langle\underline{dim}\,V,\beta\rangle=0$, where $\langle\cdot,\cdot\rangle$ is the Euler form of $Q$. Let
$$
0\longrightarrow P_1\stackrel{d^V}{\longrightarrow}
P_0\longrightarrow V\longrightarrow 0
$$
be a projective resolution of $V$. We
define the semi-invariant $c^V:=det(Hom_Q(d^V,\cdot))$ of  $ SpSI(Q,\beta)$ and $
OSI(Q,\beta)$ (see
  \cite{dw1} and \cite{s}) and, when it is possible, the semi-invariant $pf^V:=Pf(Hom_Q(d^V,\cdot))$\\
   Let $C^+$ be the Coxeter functor and let $\nabla$ be the duality functor.
    We will
  prove in the symmetric case the following
  \begin{teo}\label{p1}
Let $(Q,\sigma)$ be a symmetric quiver of tame type and let
$\beta$ be a regular symmetric dimension vector. The ring $SpSI(Q,\beta)$
is generated by semi-invariants
\begin{itemize}
\item[(i)] $c^V$ if $V\in Rep(Q)$ is such that
$\langle\underline{dim}\,V,\beta\rangle=0$,
\item[(ii)] $pf^V$ if $V\in Rep(Q)$ is such that
$\langle\underline{dim}\,V,\beta\rangle=0$, $C^+ V=\nabla V$ and
the almost split sequence $0\rightarrow\nabla V\rightarrow
Z\rightarrow V\rightarrow 0$ has the middle term $Z$ in $ORep(Q)$.
\end{itemize}
\end{teo}
\begin{teo}\label{p2}
Let $(Q,\sigma)$ be a symmetric quiver of tame type and let
$\beta$ be a regular symmetric dimension vector. The ring $OSI(Q,\beta)$
is generated by semi-invariants
\begin{itemize}
\item[(i)] $c^V$ if $V\in Rep(Q)$ is such that
$\langle\underline{dim}\,V,\beta\rangle=0$,
\item[(ii)] $pf^V$ if $V\in Rep(Q)$ is such that
$\langle\underline{dim}\,V,\beta\rangle=0$, $C^+ V=\nabla V$ and
the almost split sequence $0\rightarrow\nabla V\rightarrow
Z\rightarrow V\rightarrow 0$ has the middle term $Z$ in
$SpRep(Q)$.
\end{itemize}
\end{teo}
Differently to the results in \cite{l}, our semi-invariants $c^V$ and $pf^V$ are explicitly related to the representations of $Q$. \\
The strategy of the proofs is the following. First we set the
technique of reflection functors on the symmetric quivers. Then we
prove that we can reduce Theorems \ref{p1} and \ref{p2}, by this
technique, to particular orientations of the symmetric quivers.
Finally, we check Theorems \ref{p1} and \ref{p2} for these
orientations, first for homogeneous regular symmetric dimension vectors $ph$ and then for every other regular symmetric dimension vector.\\
In the first section we give general notions and results about
usual and symmetric quivers and their representations. We state
main results \ref{p1} and \ref{p2} and we recall some results about
representations of general linear groups and about invariant
theory.\\
In the second section we adjust to symmetric quivers the technique
of reflection functors and we prove general results about
semi-invariants of symmetric quivers.\\
In the third section we check that we can reduce Theorems \ref{p1} and
\ref{p2} to a particular orientation of symmetric quivers of
tame type. Then, using classical invariant theory and the
technique of Schur functors, we prove Theorems
\ref{p1} and \ref{p2} for symmetric quivers of tame type with
this orientation, restricting us to homogeneous regular symmetric dimension vectors. Finally, after having defined and described the symplectic and orthogonal generic decomposition of a symmetric dimension vector (for classical definition see remark 2.8 (a) in \cite{k1}), we prove  Theorems \ref{p1} and \ref{p2} for every other regular symmetric dimension vector.

\section*{Acknowledgments}
I am deeply grateful to my PhD thesis advisor Jerzy Weyman for his constant guidance and for all suggestions and ideas he shared me. I am
grateful to Elisabetta Strickland for the great
helpfulness she showed me during the three years of my graduate studies. I also
wish to thank Fabio Gavarini and Alessandro
D'Andrea for several useful discussions.

\section{Preliminary results}
Throughout the paper $\Bbbk$ denotes an algebraically closed field
of characteristic 0.
\subsection{Representations of quivers}
A quiver $Q$ is a pair $(Q_0,Q_1)$ where $Q_0$ is a set of
vertices, $Q_1$ is a set of arrows. For each arrow $a\in Q_1$ we
shall call $ta\in Q_0$ the tail of $a$ and $ha\in Q_0$ the head of
$a$. Throughout the paper we consider quivers $Q$ without oriented
cycles, i.e. in which there are no paths
$a_1\cdots a_n$ such that $ta_1=ha_n$.\\
A representation $V$ of $Q$ is a pair $\{\{V(x)\}_{x\in Q
_0},\{V(a):V(ta)\rightarrow V(ha)\}_{a\in Q_1}\}$ where $V(x)$ is
a finite dimensional vector space for every $x\in Q_0$ and $V(a)$
 is a linear map for every $a\in Q_1$. The dimension vector
of $V$ is the vector $\underline{dim}(V):=(dim V(x))_{x\in
Q_0}\in\mathbb{N}^{Q_0}$.\\ 
For a dimension vector $\alpha$ we
define
$
Rep(Q,\alpha):=\bigoplus_{a\in
Q_0}Hom(\Bbbk^{\alpha(ta)},\Bbbk^{\alpha(ha)})
$
the variety of representations of $Q$ of dimension $\alpha$.
Moreover we define the action of the group
$
SL(Q,\alpha)=\prod_{x\in Q_0}SL(\alpha(x))<\prod_{x\in Q_0}GL(\alpha(x))=GL(Q,\alpha)
$
on $Rep(Q,\alpha)$ by
$
g\cdot V=\{g_{ha}V(a)g_{ta}^{-1}\}_{a\in Q_1}
$
where $g=(g_x)_{x\in Q_0}\in GL(Q,\alpha)$ and $V\in
Rep(Q,\alpha)$.\\
A morphism $f:V\rightarrow W$ of two representations is a family
of linear maps $\{f(x):V(x)\rightarrow
W(x)\}_{x\in Q_0}$ such that $f(ha)V(a)=W(a)f(ta)$ for each $a\in Q_1$. Representations and morphisms between representations of a quiver $Q$ define respectively the objects and the morphisms of a category, denoted by $Rep(Q)$. We
denote the space of morphisms from $V$ to $W$ by $Hom_Q(V,W)$ and
the space of extensions of $V$ by $W$ by $Ext_Q^1(V,W)$.\\
Finally we define, for every $\alpha,\beta\in\mathbb{Z}^{Q_0}$,
the non symmetric bilinear form on the space $\mathbb{Z}^{Q_0}$ by
$$
\langle\alpha,\beta\rangle=\sum_{x\in
Q_0}\alpha(x)\beta(x)-\sum_{a\in Q_1}\alpha(ta)\beta(ha),
$$
called the Euler form of $Q$.\\
A vertex $x\in Q_0$ is said to be a sink (resp. a source) of $Q$
if
$x=ha$ (resp. $x=ta$) for every $a\in Q_1$ connected to $x$.\\
Let $x\in Q_0$ be a sink (resp. a source) of a quiver $Q$ and let
$\{a_1,\ldots,a_k\}$ be the arrows connected to $x$. We define the
quiver $c_x(Q)$ as follows
$$c_x(Q)_0=Q_0\quad\textrm{and}\quad c_x(Q)_1=\{c_x(a)|a\in Q_1\}$$
where $tc_x(a_i)=ha_i$,   $hc_x(a_i)=ta_i$ for every
$i\in\{1,\ldots,k\}$ and $tc_x(b)=tb$, $hc_x(b)=hb$ for every
$b\in Q_1\setminus\{a_1,\ldots,a_k\}$. Moreover, we can define the
reflection $c_x:\mathbb{Z}^{Q_0}\rightarrow\mathbb{Z}^{Q_0}$ given
for $\alpha\in\mathbb{Z}^{Q_0}$ by formula
$$
c_x(\alpha)(y)=\left\{\begin{array}{ll}
\alpha(y) & if\quad y\ne x\\
\sum_{i=1}^k \alpha(ta_i)-\alpha(x) & \textrm{otherwise}.
\end{array}\right.
$$
Finally it is known (see \cite{bgp} and \cite{dr}) that for every
$x\in Q_0$ sink or a source of $Q$ we can define respectively the
functors
$$
C_x^+:Rep(Q)\rightarrow
Rep(c_x(Q))\quad\textrm{and}\quad
C_x^-:Rep(Q)\rightarrow Rep(c_x(Q))$$ called
\textit{reflection functors}.
\begin{definizione}
Let $Q$ be a quiver with $n$ vertices without oriented cycles. We
choose the numbering $(x_1,\ldots,x_n)$ of vertices such that
$ta>ha$ for every $a\in Q_1$. We define
$$
C^+:=C^+_{x_n}\cdots C^+_{x_1}\quad\textrm{and}\quad
C^-:=C^-_{x_1}\cdots C^-_{x_n}.
$$
The functors $C^+,C^-:Rep(Q)\rightarrow Rep(Q)$ are called Coxeter
functors. \end{definizione} One proves that these functors don't
depend on the choice of numbering of vertices (see \cite{ass}
chap. VII Lemma 5.8).
\begin{teorema}[Auslander-Reiten 1975]\label{ar}
Let $Q$ be a quiver without oriented cycles.
\begin{itemize}
\item[1)] For every indecomposable non-projective representation
$V$ of $Q$ there is an almost split sequence $0\rightarrow C^+
V\rightarrow X\rightarrow V\rightarrow 0$ in $Rep(Q)$.
\item[2)] For every indecomposable non-injective
representation $V$ of $Q$ there is an almost split sequence
$0\rightarrow V\rightarrow Z\rightarrow C^-V\rightarrow 0$ in
$Rep(Q)$.
\end{itemize}
\end{teorema}
\begin{proof} See
\cite{ass}, sec. IV.3 Theorem 3.1.\end{proof}
\begin{definizione}
A quiver $Q$ is called of tame type if the underlying graph of $Q$
is of type $\widetilde{A},\widetilde{D}$ or $\widetilde{E}$.
\end{definizione}
For all next results we refer to \cite{dr}. We also consider $C^+$ and
$C^-$ as linear transformations on the space of dimension vectors,
i.e. if $V$ is a representation of a quiver with dimension
$\alpha$ then $C^{\pm}\alpha=\underline{dim}(C^{\pm}V)$.
\begin{proposizione}\label{h}
Let $Q$ be a quiver of tame type, then the quadratic form
$q_Q:\mathbb{Z}^{Q_0}\rightarrow\mathbb{Z}$ defined by
$$
q_Q(\alpha):=\sum_{x\in Q_0}\alpha(x)^2-\sum_{a\in
Q_1}\alpha(ta)\alpha(ha)
$$
is positive semi-definite and there exists a unique vector
$h\in\mathbb{N}^{Q_0}$ such that $\mathbb{Z}h$ is the radical of
$q_Q$. For quivers of type
$\widetilde{A}$ and $\widetilde{D}$ the vector $h$ has the
following form
\begin{equation}
\begin{array}{ccc}
\begin{array}{cccccccccccc}
&&1&\cdots&1&\\
\widetilde{A}:&1&&&&1,\\
&&1&\cdots&1&
\end{array}&\quad&
\begin{array}{cccccccccccc}
&1&&&&1\\
\widetilde{D}:&&2&\cdots&2&\\
&1&&&&1
\end{array}
\end{array}
\end{equation}
\end{proposizione}
\begin{definizione}
Let $V$ be an indecomposable representation of $Q$.
\begin{itemize}
\item[(i)] $V$ is preprojective if and only if $(C^+)^iV=0$ for
$i>>0$.
\item[(ii)] $V$ is preinjective if and only if $(C^-)^iV=0$ for
$i>>0$.
\item[(iii)] $V$ is regular if and only if $(C^+)^iV\neq 0$
for every $i\in\mathbb{Z}$.
\end{itemize}
\end{definizione}
\begin{definizione}
Let $V$ be a representation of $Q$. The linear map
$$
\partial:\mathbb{N}^{Q_0}\longrightarrow\mathbb{Z}
$$
defined by $\partial(\underline{dim}\,V):=\langle
h,\underline{dim}\,V\rangle$ is called defect of $V$.
\end{definizione}
\begin{lemma}\label{difetto}
Let $V$ an indecomposable representation of $Q$. $V$ is
preprojective, preinjective or regular if and only if the defect
of $V$ is respectively negative, positive or zero.
\end{lemma}
The regular representations of $Q$ form an Abelian category
$Reg(Q)$. Moreover $Reg(Q)$ is serial, i.e. every
indecomposable regular representation has only one regular
composition series and so it is only determined by its regular
socle and by its regular length.
\begin{definizione}
A simple regular module $E$ is called homogeneous if and only if
$\underline{dim}\,E=h$.
\end{definizione}
\begin{proposizione}\label{orbitetau}
Let $Q$ be a quiver of tame type. Then there exist at most three
$C^+$-orbits $\Delta=\{e_i|\,i\in I=\{0,\ldots,u\}\}$,
$\Delta'=\{e'_i|\,i\in I'=\{0,\ldots,v\}\}$,
$\Delta''=\{e''_i|\,i\in I''=\{0,\ldots,w\}\}$, of dimension
vectors of non-homogeneous simple regular representations of $Q$
($I$, $I'$, $I''$ could be empty). We can assume that
$C^+(e_i)=e_{i+1}$ for $i\in I$ ($e_{u+1}=e_0$),
$C^+(e'_i)=e'_{i+1}$ for $i\in I'$ ($e'_{v+1}=e'_0$) and
$C^+(e''_i)=e''_{i+1}$ for $i\in I''$ ($e''_{w+1}=e_0$).
\end{proposizione}
\noindent Graphycally we can present $\Delta$, $\Delta'$ and $\Delta''$ respectively as the polygons 
$$
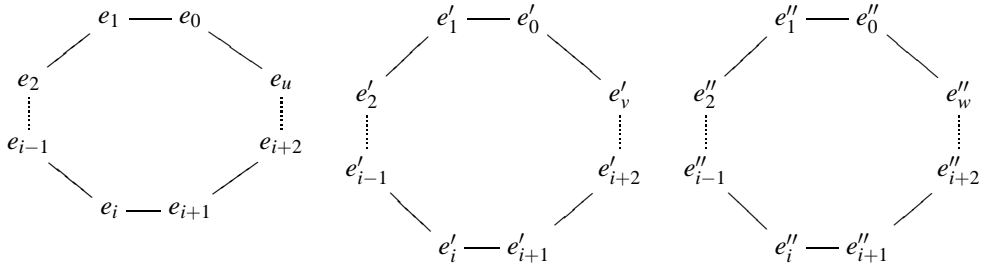
\begin{array}{ccc}
\xymatrix@-1pc{
&e_1\ar@{-}[r]\ar@{-}[dl]&e_0\ar@{-}[dr]&\\
e_2\ar@{.}[d]&&&e_u\ar@{.}[d]\\
e_{i-1}\ar@{-}[dr]&&&\ar@{-}[dl]e_{i+2}\\
&e_{i}\ar@{-}[r]&e_{i+1}&}&
\xymatrix@-1pc{
&e'_1\ar@{-}[r]\ar@{-}[dl]&e'_0\ar@{-}[dr]&\\
e'_2\ar@{.}[d]&&&e'_v\ar@{.}[d]\\
e'_{i-1}\ar@{-}[dr]&&&\ar@{-}[dl]e'_{i+2}\\
&e'_{i}\ar@{-}[r]&e'_{i+1}&}&
\xymatrix@-1pc{
&e''_1\ar@{-}[r]\ar@{-}[dl]&e''_0\ar@{-}[dr]&\\
e''_2\ar@{.}[d]&&&e''_w\ar@{.}[d]\\
e''_{i-1}\ar@{-}[dr]&&&\ar@{-}[dl]e''_{i+2}\\
&e''_{i}\ar@{-}[r]&e''_{i+1}&}
\end{array}
$$
Every dimension vector $d$ of a regular representation of $Q$ can be
decomposed uniquely as
\begin{equation}\label{dcsolita}
d=ph+\sum_{i\in I}p_ie_i+\sum_{i\in I'}p'_ie'_i+\sum_{i\in
I''}p''_ie''_i
\end{equation}
for some $p,p_i,p'_i,p''_i\in\mathbb{N}$ such that at least one of
coefficients in each family $\{p_i|\,i\in I\}$, $\{p'_i|\,i\in
I'\}$, $\{p''_i|\,i\in I''\}$ is zero.\\
We observe that the category $Reg(Q)$ can be decomposed as
direct sum of categories $\mathcal{R}_t$, with
$t=(\varphi,\psi)\in\mathbb{P}_1(\Bbbk)$. In all categories
$\mathcal{R}_t$, but at most three of these, there is only one
simple object $V_t$ which is necessarily homogeneous.
\begin{definizione}\label{indregij}
(1) We call $E_i$, $E_i'$ and $E_i''$ the simple non-homogeneous
regular representations respectively of dimension $e_i$, $e_i'$
and $e_i''$.\\
(2) We call $V_{(\varphi,\psi)}$, where $(\varphi,\psi
)\in\mathbb{P}_1(\Bbbk)$, the indecomposable regular
representation of dimension $h$.\\
(3) We define $E_{i,j}$ to be the indecomposable regular
representations with socle $E_i$ and dimension $\sum_{k=i}^j e_k$,
where $e_k$ are vertices of the arc with clockwise orientation
$\xymatrix@-1pc{e_i\ar@{-}[r]&\ar@{.}[r]&\ar@{-}[r]&e_j}$ in
$\Delta$, without repetitions of $e_k$. We denote $E_i:=E_{i,i}$
and similarly we define $E'_{i,j}$ and $E''_{i,j}$.
\end{definizione}
To simplify the notation, sometimes we will consider, without advance notice, $I=\{1,\ldots,u+1\}$ (respectively  $I'=\{1,\ldots,v+1\}$), identifying $E_{i,u+1}$ with $E_{i,0}$ and $E_{u+1,i}$ with $E_{0,i}$ for every $i\in I$ (respectively $E'_{i,v+1}$ with $E'_{i,0}$ and $E'_{v+1,i}$ with $E'_{0,i}$ for every $i\in I'$).

 \subsection{Symmetric quivers}
 \begin{definizione}
A symmetric quiver is a pair $(Q,\sigma)$ where $Q$ is a quiver
(called underlying quiver of $(Q,\sigma)$) and $\sigma$ is an
involution on $Q_0\sqcup Q_1$ such that
\begin{itemize}
   \item[(i)] $\sigma(Q_0)=Q_0$ and $\sigma(Q_1)=Q_1$,
   \item[(ii)] $t\sigma(a)=\sigma(ha)$ and $h\sigma(a)=\sigma(ta)$ for all $a\in
   Q_1$,
   \item[(iii)] $\sigma(a)=a$ whenever $a\in Q_1$ and $\sigma(ta)=ha$.
   \end{itemize}
 \end{definizione}
 Let $V$ be a representation of the underlying quiver $Q$ of a
 symmetric quiver $(Q,\sigma)$. We define the duality functor
 $\nabla:Rep(Q)\rightarrow Rep(Q)$ such that $\nabla
 V(x)=V(\sigma(x))^*$ for every $x\in Q_0$ and $\nabla
 V(a)=-V(\sigma(a))^*$ for every $a\in Q_1$. If $f:V\rightarrow W$
 is a morphism of representations $V,W\in Rep(Q)$, then $\nabla
 f:\nabla W\rightarrow \nabla V$ is defined by $\nabla
 f(x)=f(\sigma(x))^*$, for every $x\in Q_0$. We call $V$
 \textit{selfdual}
 if $\nabla V=V$.
  \begin{definizione}\label{defortsp}
 An orthogonal (resp. symplectic) representation of a symmetric quiver $(Q,\sigma)$ is a pair $(V,<\cdot,\cdot >)$, where $V$ is
  a representation of the underlying quiver $Q$ with a nondegenerate symmetric (resp. skew-symmetric) scalar product $<\cdot,\cdot >$ on $\bigoplus_{x\in Q_0}V(x)$
   such that
   \begin{itemize}
   \item[(i)] the restriction of $<\cdot,\cdot>$ to $V(x)\times V(y)$ is 0 if $y\neq\sigma(x)$,
   \item[(ii)] $<V(a)(v),w>+<v,V(\sigma(a))(w)>=0$ for all $v\in V(ta)$ and all $w\in V(\sigma(a))$.
   \end{itemize}
  \end{definizione}
  By properties \textit{(i)} and \textit{(ii)} of definition \ref{defortsp}, an orthogonal or symplectic representation\\ $(V,<\cdot,\cdot>)$ of a symmetric quiver is
  selfdual.\\
 We shall say that a dimension vector $\alpha$ is \textit{symmetric} if
 $\alpha(x)=\alpha(\sigma(x))$ for every $x\in Q_0$. Since each
 orthogonal or symplectic representation is selfdual, then
 dimension vector of an orthogonal (resp. symplectic)
 representation, which we shall call respectively orthogonal and symplectic dimension vector, is symmetric.
\begin{definizione}\label{defindspo}
An orthogonal (respectively symplectic) representation is called
indecomposable orthogonal (respectively indecomposable symplectic)
if it cannot be expressed as a direct sum of orthogonal
(respectively symplectic) representations.
\end{definizione}
 \begin{definizione}
A symmetric quiver is said to be of tame representation type if is
not of finite representation type, but in every dimension vector
the indecomposable orthogonal (symplectic) representations occur
in families of dimension $\leq 1$.
\end{definizione}
Derksen and Weyman classified the symmetric quiver of tame type in
\cite{dw2}
\begin{proposizione}\label{ctt}
Let $(Q,\sigma)$ be a symmetric tame quiver with $Q$ connected.
Then $(Q,\sigma)$ is one of the following symmetric quivers.
\begin{itemize}
\item[(1)]  Of type
$$
 \widetilde{A}^{2,0,1}_n:
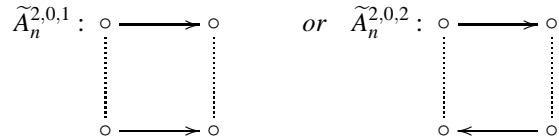
\xymatrix@-1pc{
\circ\ar[rr]\ar@{.}[dd]&&\circ\ar@{.}[dd]\\
&&&\\
 \circ\ar[rr]&&\circ}\quad \textrm{or} \quad \widetilde{A}^{2,0,2}_n:\xymatrix@-1pc{
\circ\ar[rr]\ar@{.}[dd]&&\circ\ar@{.}[dd]\\
&&&\\
 \circ&&\ar[ll]\circ}
  $$
with arbitrary orientation reversed under $\sigma$ if
$Q=\widetilde{A}_{2n+1}$ ($n\geq 1$). Here $\sigma$ is a reflection
with respect to a central vertical line (so $\sigma$ fixes two
arrows and no vertices).
\item[(2)] Of type 
$$
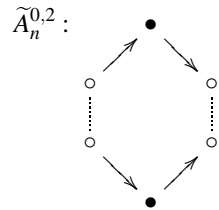
\widetilde{A}^{0,2}_n:\xymatrix@-1pc{&\bullet\ar[dr]&\\
\circ\ar[ur]\ar@{.}[d]&&\circ\ar@{.}[d]\\
\circ\ar[dr]&&\circ\\
&\bullet\ar[ur]&}
$$
with arbitrary orientation reversed under $\sigma$ if
$Q=\widetilde{A}_{2n-1}$ ($n\leq 1$). Here $\sigma$ is a
reflection with respect to a central vertical line (so $\sigma$
fixes two vertices and no arrows).
\item[(3)] Of type
$$
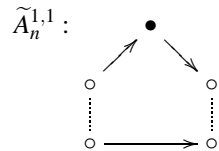
\widetilde{A}^{1,1}_n:\xymatrix@-1pc{&\bullet\ar[dr]&\\
\circ\ar[ur]\ar@{.}[d]&&\circ\ar@{.}[d]\\
\circ\ar[rr]&&\circ}
$$
with arbitrary orientation reversed under $\sigma$ if
$Q=\widetilde{A}_{2n}$ ($n\geq 1$). Here $\sigma$ is a reflection
with respect to a central vertical line (so $\sigma$ fixes one
arrow and one vertex).
\item[(4)] Of type
$$
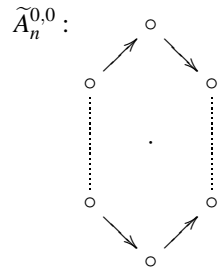
\widetilde{A}^{0,0}_n:\xymatrix@-1pc{&\circ\ar[dr]&\\
\circ\ar[ur]\ar@{.}[dd]&&\circ\ar@{.}[dd]\\
&\cdot&\\
\circ\ar[dr]&&\circ\\
&\circ\ar[ur]&}
$$
with arbitrary orientation reversed under $\sigma$ if
$Q=\widetilde{A}_{2n+1}$ ($n\geq 1$). Here $\sigma$ is a central
symmetry (so $\sigma$ fixes neither arrows nor vertices).
\item[(5)] Of type 
$$
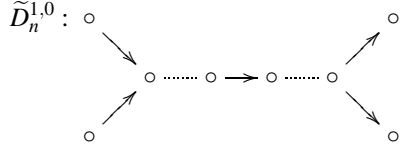
\widetilde{D}^{1,0}_n:\xymatrix@-1pc{\circ\ar[dr]&&&&&\circ\\
&\circ\ar@{.}[r]&\circ\ar[r]&\circ\ar@{.}[r]&\circ\ar[ur]\ar[dr]&\\
\circ\ar[ur]&&&&&\circ}
$$
 with arbitrary orientation reversed
under $\sigma$ if $Q=\widetilde{D}_{2n}$ ($n\geq 2$). Here
$\sigma$ is a reflection with respect to a central vertical line
(so $\sigma$ fixes one arrow and no vertices).
\item[(6)] Of type
$$
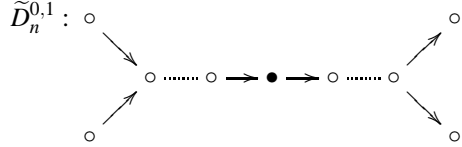
\widetilde{D}^{0,1}_n:\xymatrix@-1pc{\circ\ar[dr]&&&&&&\circ\\
&\circ\ar@{.}[r]&\circ\ar[r]&\bullet\ar[r]&\circ\ar@{.}[r]&\circ\ar[ur]\ar[dr]&\\
\circ\ar[ur]&&&&&&\circ}
$$
 with arbitrary orientation reversed
under $\sigma$ if $Q=\widetilde{D}_{2n-1}$ ($n\geq 2$). Here
$\sigma$ is a reflection with respect to a central vertical line
(so $\sigma$ fixes one vertex and no arrows).
\end{itemize}
\end{proposizione}
\begin{proof} See \cite{dw2}, proposition 4.3.\end{proof}
We describe the space of orthogonal (resp. symplectic)
representations of a symmetric quiver $(Q,\sigma)$.\\
We denote
 $ Q^{\sigma}_0$
 (respectively $Q^{\sigma}_1$) the set of vertices (respectively
arrows) fixed by $\sigma$. Thus we have partitions
$$
Q_0 = Q^+_ 0 \sqcup Q^{\sigma}_ 0 \sqcup Q^-_0
$$
$$
 Q_1 = Q^+_ 1 \sqcup Q^{\sigma}_ 1 \sqcup
Q^-_1
$$
 such that $Q^-_0 = \sigma(Q^+_ 0 )$ and $Q^-_1 = \sigma(Q^+_ 1 )$, satisfying:
 \begin{itemize}
\item[i)] $\forall a \in Q^+_ 1$ , either $\{ta, ha\} \subset Q^+_ 0$ or one of the elements in $\{ta, ha\}$ is
in $Q^+_ 0$ while the other is in $Q^{\sigma}_ 0$;
\item[ii)] $\forall x\in Q^+_0$, if $a\in Q_1$ with $ta = x$ or $ha = x$, then $a\in Q^+_ 1 \sqcup
Q^{\sigma}_1$.
\end{itemize}
 \begin{definizione}\label{delta}
 Let $(Q,\sigma)$ be a symmetric quiver. We define a linear map $\delta:\mathbb{N}^{Q_0}\rightarrow\mathbb{N}^{Q_0}$ by setting $\{\delta\alpha(i)\}_{i\in Q_0}=\{\alpha(\sigma(i))\}_{i\in Q_0}$ for every dimension vector
  $\alpha$.
   \end{definizione}
    \begin{oss}
    \begin{itemize}
\item[(i)]  Since $\sigma$ is an involution, also $\delta$ is one.
\item[(ii)] If $V$ is a representation of dimension $\alpha$ then
  $\delta\alpha=\underline{dim}(\nabla V)$. In particular if $V$
  is an orthogonal or symplectic representation of $(Q,\sigma)$ of
  dimension $\alpha$, then $\delta\alpha=\alpha$. 
\item[(iii)] If $\alpha$ and $\beta$ are
  dimension vectors, then
 \begin{eqnarray}
  \langle\alpha,\beta\rangle=\langle\delta\beta,\delta\alpha\rangle.
  \end{eqnarray}
\end{itemize}
\end{oss}

For next statements, see section 2 in \cite{dw2}.
Let $ORep(Q,\alpha)$ be the space of orthogonal $\alpha$-dimensional
  representations of a symmetric quiver $(Q,\sigma)$; by property \textit{(ii)} of definition \ref{defortsp} we have
  \begin{equation}
   ORep(Q,\alpha)\cong\bigoplus_{a\in Q^{+}_1 }Hom(\Bbbk^{\alpha(ta)},\Bbbk^{\alpha(ha)})\oplus\bigoplus_{a\in
Q^{\sigma}_1}\bigwedge^2(\Bbbk^{\alpha(ta)})^*.
  \end{equation}
  Let $SpRep(Q,\alpha)$ be the space of symplectic $\alpha$-dimensional representations of a symmetric quiver $(Q,\sigma)$; by property \textit{(ii)} of definition \ref{defortsp} we have
 \begin{equation}
 SpRep(Q,\alpha)\cong\bigoplus_{a\in Q^+_1}Hom(\Bbbk^{\alpha(ta)},\Bbbk^{\alpha(ha)})\oplus\bigoplus_{a\in Q^{\sigma}_1 }S_2(\Bbbk^{\alpha(ta)})^*.
   \end{equation}
  By property \textit{(i)} of definition \ref{defortsp}, the subgroup of $SL(Q,\alpha)$ which stabilizes $ORep(Q,\alpha)$ is
\begin{equation}
 SO(Q,\alpha)=\prod_{x\in
Q^+_0}SL(\alpha(x))\times\prod_{x\in Q^{\sigma}_0}SO(\alpha(x)),
   \end{equation}
  where $SO(\alpha(x))$ is the group of special orthogonal
  transformations for the symmetric form
  $<\cdot,\cdot>$ restricted to $V(x)$.\\
  Assuming that $\alpha(x)$ is even
  for every $x\in Q_0^{\sigma}$,  by property \textit{(i)} of definition \ref{defortsp}, the subgroup of $SL(Q,\alpha)$ which stabilizes $SpRep(Q,\alpha)$ is
\begin{equation}
   SSp(Q,\alpha)=\prod_{x\in Q^+_0}SL(\alpha(x))\times\prod_{x\in
   Q^{\sigma}_0}Sp(\alpha(x)),
   \end{equation}
where $Sp(\alpha(x))$ is the group of isometric
  transformations for the skew-symmetric form
  $<\cdot,\cdot>$ restricted to $V(x)$.\\
The action of these groups is defined by
$$
g\cdot V=\{g_{ha}V(a)g_{ta}^{-1}\}_{a\in Q_1^+\cup Q_1^{\delta}}
$$
where $g=(g_x)_{x\in Q_0}\in SO(Q,\alpha)$ (respectively
  $g\in SSp(Q,\alpha)$) and $V\in ORep(Q,\alpha)$ (respectively in $SpRep(Q,\alpha)$). In particular we can suppose $g_{\sigma(x)}=(g_x^{-1})^t$
  for every $x\in Q_0$.

\subsection{Semi-invariants of quivers without oriented cycles and main
  results}In this section first we define semi-invariants which
  appear in main results of paper and we describe some property of
  these for any quiver $Q$; then we state main theorems.\\
  Let $Q$ be a quiver with $n$ vertices. We denote
$$
SI(Q,\alpha)=\Bbbk[Rep(Q,\alpha)]^{SL(Q,\alpha)}
$$
the ring of semi-invariants of a quiver $Q$.\\
For every $g\in GL(Q,\alpha)$ the character $\tau$ at $g$ is
$\tau(g)=det(g)^{\chi_1}\cdots det(g)^{\chi_n}$, where
$\chi=(\chi_1,\ldots,\chi_n)\in \mathbb{Z}^n$ is also called
\textit{weight} if $\tau$ is a weight for some semi-invariant. So
the ring $SI(Q,\alpha)$ decomposes in graded components as
  $$
  SI(Q,\alpha)=\bigoplus_{\tau\in char(GL(Q,\alpha))}SI(Q,\alpha)_{\tau}
  $$
 where $SI(Q,\alpha)_{\tau} =\big\{f\in\Bbbk[Rep(Q,\alpha)]|g\cdot f=\tau(g)f\;\forall g\in
 GL(Q,\alpha)\big\}$.\\
For every $V\in Rep(Q,\alpha)$, we can construct a projective
resolution, called \textit{canonical resolution of} $V$:
\begin{equation}\label{Rr}
0\longrightarrow\bigoplus_{a\in Q_1}V(ta)\otimes
P_{ha}\stackrel{d^V}{\longrightarrow}\bigoplus_{x\in
Q_0}V(x)\otimes
P_{x}\stackrel{p_V}{\longrightarrow}V\longrightarrow 0
\end{equation}
where $P_x$ is the indecomposable projective associated to vertex
$x$ for every $x\in Q_0$ (see \cite{r1}), $d^V|_{V(ta)\otimes
P_{ha}}(v\otimes e_{ha})=V(a)(v)\otimes e_{ha}-v\otimes a$ and
$p_V|_{V(x)\otimes P_{x}}(v)=v\otimes e_{x}$. Applying the functor
$Hom_Q(\cdot,W)$ to $d^V$ for $W\in Rep(Q,\beta)$, we have that
the matrix associated to $Hom_Q(d^V,W)$ is square if and only if
$\langle\alpha,\beta\rangle=0$ (see \cite{s} Lemma 1.2).
\begin{definizione}
For $V\in Rep(Q,\alpha)$ such that $\langle\alpha,\beta\rangle=0$,
where $\beta\in\mathbb{N}^n$, we define
$$
\begin{array}{rcl}
c^V:Rep(Q,\beta)&\longrightarrow&\Bbbk\\
W&\longmapsto&c^V(W)=det(Hom_Q(d^V,W)).
\end{array}
$$
These are semi-invariants of weight $\langle\alpha,\cdot\rangle$,
called Schofield semi-invariants (see \cite{s} Lemma 1.4).
\end{definizione}
\begin{oss}\label{qp}
\begin{itemize}
\item[(i)] If
$$
0\rightarrow P_1\stackrel{\varphi}\rightarrow P_0\rightarrow V\rightarrow 0
$$
is another projective resolution of $V$ and $\langle\alpha,\beta\rangle=0$, than
$$
det(Hom_Q(\varphi,\cdot))=k\cdot det(Hom_Q(d^V,\cdot))
$$
for some $k\in\Bbbk$ (see \cite{s} lemma 1.2). So any projective resolution of $V$ can be used to
calculate $c^V$ (see \cite{s}). Moreover if $P$ is a projective
representation, then $c^P=0$.
\item[(ii)] If
$\langle\underline{dim}V,\underline{dim}W\rangle=0$, then
$c^V(W)=0$ if and only if $Hom_Q(V,W)\neq 0$ (see \cite{dw1}).
\end{itemize}
\end{oss}
Now we formulate the result of Derksen and Weyman about the set of
generators of the ring of semi-invariants of a quiver without
oriented cycles $Q$.
\begin{teorema}[Derksen-Weyman]\label{dw}
Let $Q$ be a quiver without oriented cycles and let $\beta$ be a
dimension vector. The ring $SI(Q,\beta)$ is spanned by
semi-invariants of the form $c^V$ of weight $\langle
\underline{dim}(V),\cdot\rangle$, for which $\langle
\underline{dim}(V),\beta\rangle=0$.
\end{teorema}
\begin{proof}
See \cite{dw1} Theorem 1.
\end{proof}
We give some property of Schofield semi-invariants.
\begin{lemma}\label{cVcV'}
Suppose that $V'$, $V$, $V''$ and $W$ are representations of $Q$,
that $\langle\underline{dim}(V),\underline{dim}(W)\rangle=0$ and
that there are exact sequences
$$
0\rightarrow V'\rightarrow V\rightarrow V''\rightarrow 0
$$
then
\begin{itemize}
\item[(i)] If $\langle\underline{dim}(V'),\underline{dim}(W)\rangle<0$, then $c^V(W)=0$
\item[(ii)] If $\langle\underline{dim}(V'),\underline{dim}(W)\rangle=0$, then $c^{V}(W)=c^{V'}(W)
c^{V''}(W)$.
\end{itemize}
\end{lemma}
\begin{proof} See \cite{dw1} Lemma 1.\end{proof}
We recall definition and properties of the
  \textit{Pfaffian} of a skew-symmetric matrix.\\
  Let $A=(a_{ij})_{1\leq i,j\leq 2n}$ be a skew-symmetric $2n\times 2n$
  matrix. Given $2n$ vectors $x_1,\ldots,x_{2n}$ in $\Bbbk^{2n}$,  we define
  $$
  F_A(x_1,\ldots,x_{2n})=\sum_{{i_1<j_1,\ldots,i_n<j_n \atop i_1<\ldots<i_n}}sgn\left(\left[\begin{smallmatrix} 1 & 2
  &\ldots & 2n-1 & 2n \\
  i_1 & j_1 & \ldots & i_n & j_n \end{smallmatrix}\right]\right)\prod_{i=1}^n(x_{s(2i-1)},x_{s(2i)}),
  $$
  where $sgn(\cdot)$
  is the sign of the permutation
  and $(\cdot,\cdot)$ is the
  skew-symmetric bilinear form associated to $A$. So $F_A$ is a
  skew-symmetric multilinear function of $x_1,\ldots,x_{2n}$.
  Since, up to a scalar, the only one skew-symmetric
  multilinear function of $2n$ vectors in $\Bbbk^{2n}$ is the
  determinant, there is a complex number $pf(A)$, called
  \textit{Pfaffian of} $A$, such that
  $$
  F_A(x_1,\ldots,x_{2n})=pf(A)det[x_1,\ldots,x_{2n}]
  $$
  where $[x_1,\ldots,x_{2n}]$ is the matrix which has the vector
  $x_i$ for $i$-th column. In particular we note that
  $$
 pf(A)=\sum_{{i_1<j_1,\ldots,i_n<j_n \atop i_1<\ldots<i_n}}sgn\left(\left[\begin{smallmatrix} 1 & 2
  &\ldots & 2n-1 & 2n \\
  i_1 & j_1 & \ldots & i_n & j_n \end{smallmatrix}\right]\right)a_{1_1j_1}\cdots
  a_{i_nj_n}.
  $$
  \begin{proposizione}
  Let $A$ be a skew-symmetric $2n\times 2n$ matrix.
  \begin{itemize}
  \item[(i)] For every invertible $2n\times 2n$ matrix $B$,
  $$
  pf(BAB^t)=det(B)pf(A);
  $$
  \item[(ii)] $det(A)=pf(A)^2$.
  \end{itemize}
  \end{proposizione}
  \begin{proof} See \cite{p}, chap. 5 sec. 3.6. \end{proof}
Let
$$
0\longrightarrow P_1\stackrel{d^V_{min}}\longrightarrow P_0\longrightarrow V\longrightarrow 0
$$ 
be the minimal projective resolution of $V\in Rep(Q,\alpha)$ and let $\beta$ be a dimension vector
such that $\langle\alpha,\beta\rangle=0$, we will prove (Lemma \ref{ss1}) that, under some hypothesis, 
$Hom_Q(d^V_{min},\cdot)$ is skew-symmetric on $Rep(Q,\beta)$, so we
can define
$$
\begin{array}{rcl}
pf^V:Rep(Q,\beta)&\longrightarrow&\Bbbk\\
W&\longmapsto&pf^V(W)=Pf(Hom_Q(d^V_{min},W)).
\end{array}
$$
In this work we describe a set of generators of the rings of
semi-invariants of symmetric quivers of tame type.\\
Let $\alpha$ be a dimension vector of an orthogonal or symplectic
representation, we denote
$$
OSI(Q,\alpha):=\Bbbk[ORep(Q,\alpha)]^{SO(Q,\alpha)}\quad\textrm{and}\quad
SpSI(Q,\alpha):=\Bbbk[SpRep(Q,\alpha)]^{SSp(Q,\alpha)}
$$
respectively the ring of orthogonal semi-invariants and the ring
of symplectic semi-invariants of a symmetric quiver
$(Q,\sigma)$.\\
We state the main theorems
\begin{teorema}\label{tp1}
Let $(Q,\sigma)$ be a symmetric quiver of tame type and let
$\beta$ be a regular symmetric dimension vector. The ring $SpSI(Q,\beta)$
is generated by semi-invariants
\begin{itemize}
\item[(i)] $c^V$ if $V\in Rep(Q)$ is such that
$\langle\underline{dim}\,V,\beta\rangle=0$,
\item[(ii)] $pf^V$ if $V\in Rep(Q)$ is such that
$\langle\underline{dim}\,V,\beta\rangle=0$, $C^+ V=\nabla V$ and
the almost split sequence $0\rightarrow\nabla V\rightarrow
Z\rightarrow V\rightarrow 0$ has the middle term $Z$ in $ORep(Q)$.
\end{itemize}
\end{teorema}
\begin{teorema}\label{tp2}
Let $(Q,\sigma)$ be a symmetric quiver of tame type let $\beta$ be
a regular symmetric dimension vector. The ring $OSI(Q,\beta)$ is generated
by semi-invariants
\begin{itemize}
\item[(i)] $c^V$ if $V\in Rep(Q)$ is such that
$\langle\underline{dim}\,V,\beta\rangle=0$,
\item[(ii)] $pf^V$ if $V\in Rep(Q)$ is such that
$\langle\underline{dim}\,V,\beta\rangle=0$, $C^+ V=\nabla V$ and
the almost split sequence $0\rightarrow\nabla V\rightarrow
Z\rightarrow V\rightarrow 0$ has the middle term $Z$ in
$SpRep(Q)$.
\end{itemize}
\end{teorema}
In section 3.1, we prove theorems \ref{tp1} and \ref{tp2} for
$\beta=ph$; in section 3.2, for any regular symmetric dimension
vector $\beta$.
\begin{definizione}\label{popssp}
Let $(Q,\sigma)$ be a symmetric quiver. We will say that $V\in
Rep(Q)$ satisfies property \textit{(Op)} if
\begin{itemize}
\item[(i)] $V=C^-\nabla V$
\item[(ii)] the almost split sequence $0\rightarrow\nabla V\rightarrow
Z\rightarrow V\rightarrow 0$ has the middle term $Z$ in $ORep(Q)$.
\end{itemize}
Similarly we will say that $V\in Rep(Q)$ satisfies property
\textit{(Spp)} if
\begin{itemize}
\item[(i)] $V=C^-\nabla V$
\item[(ii)] the almost split sequence $0\rightarrow\nabla V\rightarrow
Z\rightarrow V\rightarrow 0$ has the middle term $Z$ in
$SpRep(Q)$.
\end{itemize}
\end{definizione}

\subsection{Invariant theory and Schur modules}
Let $G$ be an algebraic group, $V$ a rational representation of
$G$ and $\Bbbk[V]$ the algebra of regular functions of $V$.\\
If $\mathcal{X}(G)$ is the set of characters of $G$, then the ring
of the semi-invariants of $G$ on $V$ is defined by
$$
SI(G,V)=\bigoplus_{\chi\in\mathcal{X}(G)}SI(G,V)_\chi
$$
where $SI(G,V)_\chi=\{f\in\Bbbk[V]|g\cdot f=\chi(g)f,\,\forall
g\in G\}$ is called weight space of weight $\chi$. The following
lemma describes $SI(G,V)$ in the case when $G$ has an open orbit
on $V$.
\begin{lemma}[Sato-Kimura]\label{sk}
Let $G$ be a connected linear algebraic group and $V$ a rational
representation of $G$. We suppose that the action of $G$ on $V$
has an open orbit. Then $SI(G,V)$ is a polynomial $\Bbbk$-algebra
and the weights of the generators of $SI(G,V)$ are linearly
independent in $\mathcal{X}(G)$. Moreover, the dimensions of the
spaces $SI(G,V)_\chi$ are 0 or 1.
\end{lemma}
\begin{proof}
See \cite{sk}, sect. 4, Lemma 4 and Proposition 5.
\end{proof}
Let $G=GL_n(\Bbbk)$ be the general linear group over $\Bbbk$.
There exists an isomorphism $\mathbb{Z}\cong \mathcal{X}(G)$ which
sends an element $a$ of $\mathbb{Z}$ in $(det)^a$ (where $det$
associates to $g\in G$ its determinant). We identify $G$ with the
group $GL(V)$ of linear automorphisms of a vector space $V$ of
dimension $n$. So we have
 $$
 SI(G,V)=\Bbbk[V]^{SL(V)}.
 $$
Let $T$ and $\mathcal{X}(T)$ respectively be the maximal torus in
$G$ (i.e. the group of diagonal matrices) and the set of
characters of $T$. The irreducible rational representations of $G$
are parametrized by the set
$$
\mathcal{X}^+(T)=\{\lambda=(\lambda_1,\ldots,\lambda_n)\in\mathbb{Z}^n|\lambda_1\geq\cdots\geq\lambda_n\}
$$
of the integral dominant weights for $GL_n(\Bbbk)$. The
irreducible rational representations $S_{\lambda}V$ of
$G=GL_n(\Bbbk)$ corresponding to the dominant weight
$\lambda\in\mathcal{X}^+(T)$ are called Schur modules. In the case
when $\lambda_n\geq 0$ (i.e. $\lambda$ is a partition of
$\lambda_1+\cdots+\lambda_n$), a description of $S_{\lambda}V$ is
given in \cite{abw} and \cite{p}. For every
$\lambda\in\mathcal{X}^+(T)$, we can define $S_{\lambda}V$ as
follows
$$
S_{(\lambda_1,\ldots,\lambda_n)}V=S_{(\lambda_1-\lambda_n,\ldots,\lambda_{n-1}-\lambda_n,0)}V\otimes(\bigwedge^n
V)^{\otimes\lambda_n}.
$$
Let $\lambda=(\lambda_1,\ldots,\lambda_n)$ be a partition of
$|\lambda|:=\lambda_1+\cdots+\lambda_n$. We call \textit{height}
of $\lambda$, denoted by $ht(\lambda)$, the number $k$ of nonzero
components of $\lambda$ and we denote the transpose of $\lambda$
by $\lambda'$.
\begin{teorema}[Properties of Schur modules]\label{pms}
Let $V$ be vector space of dimension $n$ and $\lambda$ be an
integral dominant weight.
\begin{itemize}
\item[(i)] $S_{\lambda}V=0\Leftrightarrow ht(\lambda)>0$.
\item[(ii)] $dim\,S_{\lambda}V=1\Leftrightarrow
\lambda=(\overbrace{k,\ldots,k}^n)$ for some $n\in\mathbb{Z}$.
\item[(iii)] $\left(S_{(\lambda_1,\ldots,\lambda_n)}V\right)^*\cong
S_{(\lambda_1,\ldots,\lambda_n)}V^*\cong
S_{(-\lambda_n,\ldots,-\lambda_1)}V.$
\end{itemize}
\end{teorema}
\begin{proof}
See Theorem 6.3 in \cite{fh}.
\end{proof}
\begin{teorema}[Cauchy formulas]\label{fc}
Let $V$ and $W$ be two finite dimensional vector spaces. Then
\begin{itemize}
\item[(i)] As representations of $GL(V)\times GL(W)$,
$$
\begin{array}{ccc}
S_d(V\otimes W)=\bigoplus_{|\lambda|=d}S_{\lambda}V\otimes
S_{\lambda}W&\quad\textrm{and}\quad& \bigwedge^d(V\otimes
W)=\bigoplus_{|\lambda|=d}S_{\lambda}V\otimes S_{\lambda'}W.
\end{array}
$$
\item[(ii)] As representations of
$GL(V)$,
$$\begin{array}{ccc}S_d(S_2(V))=\bigoplus_{|\lambda|=d}S_{2\lambda}V&\quad\textrm{and}\quad& S_d(\bigwedge^2(V))
=\bigoplus_{|\lambda|=d}S_{2\lambda'}V,
\end{array}
$$
where $2\lambda=(2\lambda_1,\ldots,2\lambda_k)$ if
$\lambda=(\lambda_1,\ldots,\lambda_k)$.
\end{itemize}
\end{teorema}
\begin{proof}
See \cite{p} chap. 9 sec. 6.3 and sec 8.4, chap 11 sec. 4.5.
\end{proof}
The decomposition of tensor product of Schur modules is
$$
S_{\lambda}V\otimes S_{\mu}V=\bigoplus_{\nu}c^{\nu}_{\lambda\mu}
S_{\nu}V,
$$
where the coefficients $c^{\nu}_{\lambda\mu}$ are called
Littlewood-Richardson coefficients.\\
There exists a combinatorial formula to calculate
$c^{\nu}_{\lambda\mu}$, called \textit{Littlewood-Richardson rule}
(see in \cite{p} chap. 12 sec. 5.3).\\
 Finally we state other two results on Schur modules and invariant theory.
 \begin{proposizione}\label{i1}
 Let $V$ be a vector space of dimension $n$.
 $$
 \begin{array}{ccc}
 (S_{\lambda}V)^{SL(V)}\neq 0&\Longleftrightarrow &\lambda=(k^n)
 \end{array}
 $$
 for some $k$ and in this case $S_{\lambda}V$, and so also $(S_{\lambda}V)^{SL(V)}$, have dimension one.
 \end{proposizione}
 \begin{proof}
See Corollary p. 388 in \cite{p}.
\end{proof}
 \begin{proposizione}\label{i2}
  Let $V$ be a vector space of dimension $n$ and let $\lambda$ and $\mu$ be two integral dominant weights. Then
  $$
  \begin{array}{ccc}
  S_{\lambda}V\otimes S_{\mu}V\neq 0
  &
  \Longleftrightarrow
 &
 \lambda_i-\lambda_{i+1}=\mu_{n-i}-\mu_{n-i+1}
 \end{array}
  $$
  for every $i\in\{1,\ldots,n-1\}$ and in this case the semi-invariant is unique (up to a non zero scalar) and has weight
  $\lambda_1+\mu_n=\lambda_2+\mu_{n-1}=\cdots=\lambda_n+\mu_1$.
  \end{proposizione}

  \begin{proof}It is a Corollary of (7.11) in \cite{m} chap. 1 sec. 5.\end{proof}

  \begin{proposizione}\label{i3}
  Let $V$ be an orthogonal space of dimension $n$ and let $W$ be a
  symplectic space of dimension $2n$.
 $$
 \begin{array}{ccc}
 dim\,(S_{\lambda}V)^{SO(V)}=\left\{\begin{array}{ll}1 &
  \textrm{if}\quad \lambda=2\mu+(k^n)\\
  0 & \textrm{otherwise}\end{array}\right. &\quad\textrm{and}\quad&
   dim\,(S_{\lambda}W)^{Sp(W)}=\left\{\begin{array}{ll}1 &
  \textrm{if}\quad \lambda=2\mu'\\
  0 & \textrm{otherwise}\end{array}\right.
  \end{array}
  $$
  for some partition $\mu$ and for some $k\in\mathbb{N}$.
  \end{proposizione}
  \begin{proof} See \cite{p} chap. 11 corollaries 5.2.1 and 5.2.2.
  \end{proof}

\section{Reflection functors and semi-invariants of symmetric quivers}
\subsection{Reflection functors for symmetric quivers}
We adjust the technique of reflection functors to symmetric
quivers. See \cite{ar} section 2.1 for the proofs of results of
this paragraph.
\begin{definizione}
Let $(Q,\sigma)$ be a symmetric quiver. A sink (resp. a source)
$x\in Q_0$ is called admissible if there are no arrows connecting
$x$ and $\sigma(x)$.
\end{definizione}
By definition of $\sigma$, $x$ is an admissible sink (resp. a source) if and
only if $\sigma(x)$ is an admissible source (resp. a sink). We call
$(x,\sigma(x))$ \textit{the admissible sink-source pair}. So we
can define $c_{(x,\sigma(x))}:=c_{\sigma(x)}c_x$. Moreover we can prove that $(c_{(x,\sigma(x))}Q,\sigma)$ is a symmetric quiver (see Lemma 2.2 in \cite{ar}).

\begin{definizione}
Let $(Q,\sigma)$ be a symmetric quiver. A sequence
$x_1,\ldots,x_m\in Q_0$ is an admissible sequence of sinks (or
sources) for admissible sink-source pairs if $x_{i+1}$ is an
admissible sink (resp. source) in $c_{(x_i,\sigma(x_i))}\cdots
c_{(x_1,\sigma(x_1))}(Q)$ for $i=1,\ldots,m-1$.
\end{definizione}

The underlying graph of $\widetilde{D}$ is a tree, so by Proposition 2.4 in \cite{ar}, applying a composition of reflections at admissible
sink-source pairs, from any orientation of $\widetilde{D}$ we can
get the following orientation
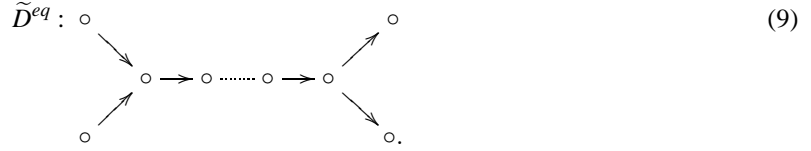
\begin{equation}\label{Dequi}
\widetilde{D}^{eq}:\xymatrix @-1pc{\circ\ar[dr] & & & & & \circ\\
& \circ\ar[r]& \circ\ar@{.}[r]&\circ\ar[r]&\circ\ar[ur]\ar[dr] &\\
\circ\ar[ur] & & & & & \circ.}
\end{equation}
\begin{definizione}\label{tipoA}
We will say that a symmetric quiver is of type $(s,t,k,l)$ if
\begin{itemize}
\item[(i)] it is of type $\widetilde{A}$,
\item[(ii)] $|Q_1^{\sigma}|=s$ and $|Q_0^{\sigma}|=t$,
\item[(iii)] it has $k$ counterclockwise arrows and $l$ clockwise arrows in $Q_1^+\sqcup
Q_1^-$.
\end{itemize}
\end{definizione}
By Proposition \ref{ctt}, $s,t\in\{0,1,2\}$ and if either $s$ or
$t$ are not zero, then $s+t=2$. Moreover, by symmetry, we note
that $k$ and $l$ have to be even.\\
One proves, restricting to subquivers of type $A$, by a simple
combinatorial argument the following (for further details, see Proposition 1.3.8 in \cite{a1})
\begin{proposizione}\label{oA}
Let $(Q,\sigma)$ be a symmetric quiver of type $\widetilde{A}$
such that $Q$ is without oriented cycles. Then there is an
admissible sequence of sinks $x_1,\ldots,x_s$ of $Q$ for
admissible sink-source pairs such that
$c_{(x_1,\sigma(x_1))}\cdots c_{(x_s,\sigma(x_s))}Q$ is one of the
quivers:
\begin{itemize}
\item[(1)]
$$
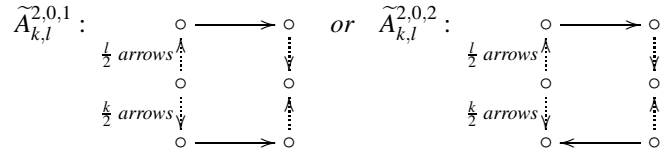
\widetilde{A}^{2,0,1}_{k,l}:\xymatrix@-1pc{
\circ\ar[rr]&&\circ\ar@{.>}[d]\\
\circ\ar@{.>}[d]_{\frac{k}{2}\;\textrm{arrows}}\ar@{.>}[u]^{\frac{l}{2}\;\textrm{arrows}}&&\circ\\
 \circ\ar[rr]&&\circ\ar@{.>}[u],}\quad\textrm{or}\quad\widetilde{A}^{2,0,2}_{k,l}:\xymatrix@-1pc{
\circ\ar[rr]&&\circ\ar@{.>}[d]\\
\circ\ar@{.>}[d]_{\frac{k}{2}\;\textrm{arrows}}\ar@{.>}[u]^{\frac{l}{2}\;\textrm{arrows}}&&\circ\\
 \circ&&\ar[ll]\circ\ar@{.>}[u],}$$
if $(Q,\sigma)$ is of type $(2,0,k,l)$;
\item[(2)]
$$
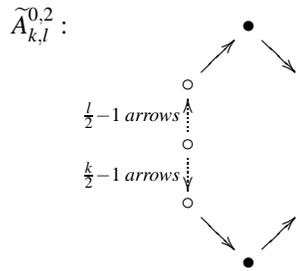
\widetilde{A}^{0,2}_{k,l}:\xymatrix@-1pc{&\bullet\ar[dr]&\\
\circ\ar[ur]&&\circ\ar@{.>}[d]\\
\circ\ar@{.>}[d]_{\frac{k}{2}-1\;\textrm{arrows}}\ar@{.>}[u]^{\frac{l}{2}-1\;\textrm{arrows}}&&\circ\\
\circ\ar[dr]&&\circ\ar@{.>}[u]\\
&\bullet\ar[ur]&,}
$$
if $(Q,\sigma)$ is of type $(0,2,k,l)$;
\item[(3)]
$$
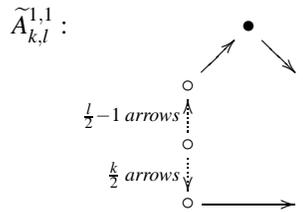
\widetilde{A}^{1,1}_{k,l}:\xymatrix@-1pc{&\bullet\ar[dr]&\\
\circ\ar[ur]&&\circ\ar@{.>}[d]\\
\circ\ar@{.>}[u]^{\frac{l}{2}-1\;\textrm{arrows}}\ar@{.>}[d]_{\frac{k}{2}\;\textrm{arrows}}&&\circ\\
 \circ\ar[rr]&&\circ\ar@{.>}[u],}
$$
if $(Q,\sigma)$ is of type $(1,1,k,l)$;
\item[(4)]
$$
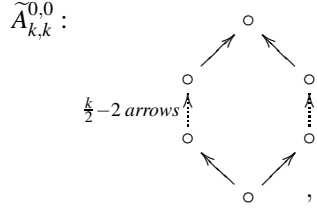
\widetilde{A}^{0,0}_{k,k}:\xymatrix@-1pc{&\circ&\\
\circ\ar[ur]&&\circ\ar[ul]\\
\circ\ar@{.>}[u]^{\frac{k}{2}-2\;\textrm{arrows}}&&\circ\ar@{.>}[u]\\
&\circ\ar[ul]\ar[ur]&,}
$$
if $(Q,\sigma)$ if of type $(0,0,k,k)$.
\end{itemize}
\end{proposizione}

Let $(Q,\sigma)$ be a symmetric quiver and $(x,\sigma(x))$ a
sink-source admissible pair. For every $V\in Rep(Q)$, we define
the reflection functors
$$
C^+_{(x,\sigma(x))}V:=C^-_{\sigma(x)}C^+_x V\quad\textrm{and}\quad
C^-_{(x,\sigma(x))}V:=C^-_xC^+_{\sigma(x)} V.
$$
\begin{proposizione}\label{C+taunabla}
Let $(Q,\sigma)$ and $(Q',\sigma)$ be two symmetric quivers with
the same underlying graph. We suppose that
$Q'=c_{(x_m,\sigma(x_m))}\cdots c_{(x_1,\sigma(x_1))}(Q)$ for some
admissible sequence of sinks $x_1,\ldots,x_m\in Q_0$ for
admissible sink-source pairs and let
$V'=C^+_{(x_m,\sigma(x_m))}\cdots C^+_{(x_1,\sigma(x_1))}V\in
Rep(Q')$. Then
$$
V=C^-\nabla V\Leftrightarrow V'=C^-\nabla V'.
$$
\end{proposizione}
\begin{proof} See Corollary 2.6 in \cite{ar}.
\end{proof}
\begin{proposizione}\label{sptosp}
Let $(Q,\sigma)$ be a symmetric quiver and let $x$ be an
admissible sink. Then $V$ is an orthogonal (resp. symplectic)
representation of $(Q,\sigma)$ if and only if
$C^+_{(x,\sigma(x))}V$ is an orthogonal (resp. symplectic)
representation of $(c_{(x,\sigma(x))}Q,\sigma)$. Similarly for
$C^-_{(x,\sigma(x))}$ if $x$ is an admissible source.
\end{proposizione}
\begin{proof} See Proposition 2.7 in \cite{ar}.
\end{proof}

\subsection{Orthogonal and symplectic semi-invariants}
If $W$ is a vector space of dimension $n$, we denote
$\widetilde{Gr}(r,W)$ the set of all decomposable tensors
$w_1\wedge\ldots\wedge w_r$, with $w_1,\ldots,w_r\in W$, inside
$\bigwedge^r W$.
\begin{lemma}\label{kac}
If $x$ is an admissible sink or source for a symmetric quiver
$(Q,\sigma)$ and $\alpha$ is a dimension vector such that
$c_{(x,\sigma(x))}\alpha(x)\geq 0$, then
\begin{itemize}
\item[i)] if $c_{(x,\sigma(x))}\alpha(x)> 0$, then there exist isomorphisms
$$
\begin{array}{ccc}
SpSI(Q,\alpha)\stackrel{\varphi^{Sp}_{x,\alpha}}{\longrightarrow}
SpSI(c_{(x,\sigma(x))}Q,c_{(x,\sigma(x))}\alpha)
&
\quad\textrm{and}\quad
&
OSI(Q,\alpha)\stackrel{\varphi^{O}_{x,\alpha}}{\longrightarrow}
OSI(c_{(x,\sigma(x))}Q,c_{(x,\sigma(x))}\alpha);
\end{array}
$$
\item[ii)] if $c_{(x,\sigma(x))}\alpha(x)= 0$, then there exist isomorphisms
$$
\begin{array}{ccc}
SpSI(Q,\alpha)\stackrel{\varphi^{Sp}_{x,\alpha}}{\longrightarrow}
SpSI(c_{(x,\sigma(x))}Q,c_{(x,\sigma(x))}\alpha)[y]
&
\quad\textrm{and}\quad
&
OSI(Q,\alpha)\stackrel{\varphi^{O}_{x,\alpha}}{\longrightarrow}
OSI(c_{(x,\sigma(x))}Q,c_{(x,\sigma(x))}\alpha)[y]
\end{array}
$$
\end{itemize}
where $R[y]$ denotes a polynomial ring with coefficients in $R$.
\end{lemma}
\begin{proof} See \cite{ar}, Lemma 2.8
\end{proof}
In next section we will see that for tame type $\varphi_{x,\alpha}^{Sp}$ and
$\varphi^{O}_{x,\alpha}$ send $c^V$ to $c^{C^+_{(x,\sigma(x))}V}$ (propositions \ref{SIcx1} and \ref{SIcx2}).\\
We recall that, by definition, symplectic groups or orthogonal
groups act
  on the spaces which are defined on the
  vertices in $Q_0^{\sigma}$, so we have
  \begin{definizione}\label{wsQ}
  Let $V$ be a
  representation of the underlying quiver $Q$ with $\underline{dim}V=\alpha$ such that
  $\langle\alpha,\beta\rangle=0$ for some symmetric
  dimension vector $\beta$. The weight of $c^V$ on
  $SpRep(Q,\beta)$ (respectively on $ORep(Q,\beta)$) is
  $\langle\alpha,\cdot\rangle-\sum_{x\in
  Q_0^{\sigma}}\varepsilon_{x,\alpha}$, where
  \begin{equation}\label{epsilon}
  \varepsilon_{x,\alpha}(y)=\left\{\begin{array}{ll}
  \langle\alpha,\cdot\rangle(x) & y=x\\
  0 & \textrm{otherwise}.\end{array}\right.
  \end{equation}
  \end{definizione}
  We shall say that a weight is \textit{symmetric} if
  $\chi(i)=-\chi(\sigma(i))$ for every $i\in Q_0$.
  \begin{oss}\label{pesi}
  Let $(Q,\sigma)$ be a symmetric quiver and $V\in Rep(Q,\alpha)$. We note that
  $$
  \langle\underline{dim}(C^-\nabla
  V),\cdot\rangle(i)=-\langle\alpha,\cdot\rangle(\sigma(i))
  $$
  for every $i\in Q_0$. So, if $C^-\nabla V=V$ then
  $\chi=\langle\alpha,\cdot\rangle$ is a symmetric weight.
  \end{oss}
\begin{lemma}\label{cV=cVnabla}
Let $(Q,\sigma)$ be a symmetric quiver. For every representation
  $V$ of the underlying quiver $Q$ and for every orthogonal or
  symplectic representation $W$ such that
  $\langle\underline{dim}(V),\underline{dim}(W)\rangle=0$, we have
  $$
  c^{V}(W)=c^{C^-\nabla V}(W).
  $$
\end{lemma}
\begin{proof} See \cite{ar}, Corollary 2.16. \end{proof}
We conclude this section with two lemma which will be useful
later.
\begin{lemma}\label{cl}
Let
$$
\xymatrix@-1pc{&\ar@{.}[dr]&&&&\ar@{.}[dr]
&&&&\\
(Q,\sigma):&\ar@{.}[r]&y\ar[r]^{a}&x\ar[r]^{b}&z\ar@{.}[r]\ar@{.}[ur]\ar@{.}[dr]&
\ar@{.}[r]&\sigma(z)\ar[r]^{a}&\sigma(x)\ar[r]^{b}&\sigma(y)\ar@{.}[r]\ar@{.}[ur]\ar@{.}[dr]&\\
&\ar@{.}[ur]&&&&\ar@{.}[ur] &&&&}
$$
be a symmetric quiver. Assume there exist only two arrows in
$Q_1^+$ incident to $x\in Q_0^+$, $a:y\rightarrow x$ and
$b:x\rightarrow z$ with $y,z\in Q_0^+\cup Q_0^{\sigma}$. Let $V$
be an orthogonal or symplectic representation with symmetric dimension vector $(\alpha_i)_{i\in Q_0}=\alpha$ such that $\alpha_x\geq max\{\alpha_y,\alpha_z\}$.\\
We define the symmetric quiver $Q'=((Q_0',Q_1'),\sigma)$ with
$n-2$ vertices such that $Q_0'=Q_0\setminus\{x,\sigma(x)\}$ and
$Q_1'=Q_1\setminus\{a,b,\sigma(a),\sigma(b)\}\cup\{ba,\sigma(a)\sigma(b)\}$
and let $\alpha'$ be the dimension of $V$ restricted to $Q'$. \\
We have:
\begin{itemize}
\item[(Sp)] Assume $V$ symplectic. Then
\begin{itemize}
\item[(a)] if $\alpha_x> max\{\alpha_y,\alpha_z\}$ then
$SpSI(Q,\alpha)=SpSI(Q',\alpha')$,
\item[(b)] if $\alpha_x=\alpha_y>\alpha_z$ then
$SpSI(Q,\alpha)=SpSI(Q',\alpha')[detV(a)]$,
\item[(b')] if $\alpha_x=\alpha_z>\alpha_y$ then
$SpSI(Q,\alpha)=SpSI(Q',\alpha')[detV(b)]$,
\item[(c)] if $\alpha_x=\alpha_y=\alpha_z$ then
 $SpSI(Q,\alpha)=SpSI(Q',\alpha')[detV(a),detV(b)]$.
\end{itemize}
\item[(O)] Assume $V$ orthogonal. Then
\begin{itemize}
\item[(a)] if $\alpha_x> max\{\alpha_y,\alpha_z\}$ then
$OSI(Q,\alpha)=OSI(Q',\alpha')$,
\item[(b)] if $\alpha_x=\alpha_y>\alpha_z$ then
$OSI(Q,\alpha)=OSI(Q',\alpha')[detV(a)]$,
\item[(b')] if $\alpha_x=\alpha_z>\alpha_y$ then
$OSI(Q,\alpha)=OSI(Q',\alpha')[detV(b)]$,
\item[(c)] if $\alpha_x=\alpha_y=\alpha_z$ then
 $OSI(Q,\alpha)=OSI(Q',\alpha')[detV(a),detV(b)]$.
\end{itemize}
\end{itemize}
\end{lemma}
\begin{proof} See \cite{ar}, Lemma 2.17.\end{proof}
Similarly one proves the following
\begin{lemma}\label{cls}
Let
$$
\xymatrix@-1pc{&\ar@{.}[dr]&&&&&\\
(Q,\sigma):&\ar@{.}[r]&y\ar[r]^{a}&x\ar[r]^{b}&\sigma(x)\ar[r]^{\sigma(a)}&\sigma(y)\ar@{.}[r]\ar@{.}[ur]\ar@{.}[dr]&\\
&\ar@{.}[ur]&&&&&}
$$
be a symmetric quiver with $n$ vertices such that there exist only
two arrows $a$ and $b$ incident to the vertex $x$ in $Q_0$ and $b$
is fixed by $\sigma$. Let $V$ be an orthogonal or symplectic
representation of $(Q,\sigma)$ with $\underline{dim}(V)=\alpha$ such
that $\alpha_x\geq\alpha_y$. Moreover we define the symmetric
quiver $(Q',\sigma)=((Q_0',Q_1'),\sigma)$ with $n-2$ vertices such
that $Q_0'=Q_0\setminus\{x,\sigma(x)\}$ and
$Q_1'=Q_1\setminus\{a,b,\sigma(a)\}\cup\{\sigma(a)ba\}$.\\
Let $\alpha'$ be the dimension of $V$ restricted to $Q'$.
\begin{itemize}
\item[(Sp)] If $V$ is symplectic, then
\begin{itemize}
\item[(i)]  $\alpha_x>\alpha_y\Longrightarrow
SpSI(Q,\alpha)=SpSI(Q',\alpha')[detV(b)]$
\item[(ii)]  $\alpha_x=\alpha_y\Longrightarrow
SpSI(Q,\alpha)=SpSI(Q',\alpha')[detV(a)]$.
\end{itemize}
\item[(O)] If $V$ is orthogonal, then
\begin{itemize}
\item[(i)]  $\alpha_x>\alpha_y\;\textrm{and}\;\alpha_x\;\textrm{is even}\Longrightarrow
OSI(Q,\alpha)=OSI(Q',\alpha')[pfV(b)]$
\item[(ii)]  $\alpha_x=\alpha_y\Longrightarrow
OSI(Q,\alpha)=OSI(Q',\alpha')[detV(a)]$.
\end{itemize}
\end{itemize}
\end{lemma}

\section{Semi-invariants of symmetric quivers of tame type}
In this section we prove theorems \ref{tp1} and \ref{tp2} for the
symmetric quivers of tame type. We recall that the underlying
quiver of a symmetric quiver of tame type is either
$\widetilde{A}$ or $\widetilde{D}$ as in Proposition \ref{ctt}. As
done for the finite case in \cite{ar}, we again reduce the
proof to particular orientations (orientations in Proposition
\ref{oA} for $\widetilde{A}$ and orientation of
$\widetilde{D}^{eq}$ for
$\widetilde{D}$).\\
We start with some result which will be useful later.
\begin{lemma}\label{ss1}
  Let $(Q,\sigma)$ be a symmetric quiver of tame type. Let $d_{min}^V$
  be the matrix of the minimal projective presentation of $V\in Rep(Q,\alpha)$ and
  let $\beta$ be a symmetric dimension vector such that
  $\langle\alpha,\beta\rangle=0$. Then
  \begin{itemize}
  \item[(1)] $Hom_Q(d_{min}^V,\cdot)$ is skew-symmetric on $SpRep(Q,\beta)$  if and only if $V$ satisfies property \textit{(Op)};
  \item[(2)] $Hom_Q(d_{min}^V,\cdot)$ is skew-symmetric on $ORep(Q,\beta)$ if and only if $V$ satisfies property \textit{(Spp)}.
  \end{itemize}
\end{lemma}
\begin{proof}
First we note, by Auslander-Reiten quiver of $Q$, that if
$(Q,\sigma)$ is a symmetric quiver of tame type, then the only
representations $V\in Rep(Q)$ such that
$C^-\nabla V=V$ are regular ones.\\
We prove $(1)$ only for symmetric quivers $(Q,\sigma)$ of type
$(1,1,k,l)$ (see Definition \ref{tipoA}), because for $(2)$ and
the other cases one proceeds similarly (for the details for symmetric quiver of type $\widetilde{D}^{0,1}_n$, see Lemma 1.4.6 in \cite{a1}).\\
We call $(Q',\sigma)$ the symmetric quiver with the same
underlying graph of $(Q,\sigma)$ and with orientation as in
Proposition \ref{oA}. By Proposition \ref{oA}, there exists a
sequence $x_1,\ldots,x_m$ of admissible sink for admissible
sink-source pairs such that $c_{(x_m,\sigma(x_m))}\cdots
c_{(x_1,\sigma(x_1))}Q=Q'$. We call
$V':=C^+_{(x_m,\sigma(x_m))}\cdots C^+_{(x_1,\sigma(x_1))}V$ for
every $V\in Rep(Q)$ and if $\alpha=\underline{dim}\,V$, then
$\alpha':=c_{(x_m,\sigma(x_m))}\cdots
c_{(x_1,\sigma(x_1))}\alpha$. We note that, by Proposition
\ref{C+taunabla} and Proposition \ref{sptosp}, $V$ satisfies
property \textit{(Op)} (respectively property \textit{(Spp)}) if
and only if $V'$ satisfies property \textit{(Op)} (respectively
property \textit{(Spp)}).\\
We consider the following labelling for
$Q'=\widetilde{A}^{1,1}_{k,l}$:
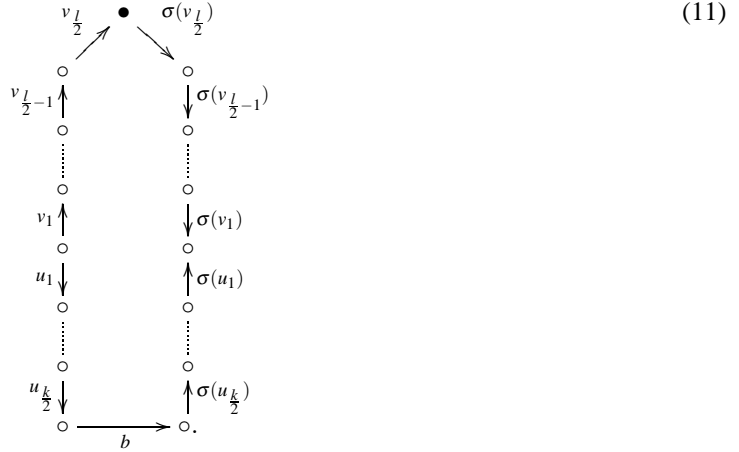
\begin{equation}\label{A11}
\xymatrix@-1pc{ &\bullet\ar[dr]^{\sigma(v_{\frac{l}{2}})}&\\
\circ\ar[ur]^{v_{\frac{l}{2}}}&&\circ\ar[d]^{\sigma(v_{\frac{l}{2}-1})}\\
\circ\ar[u]^{v_{\frac{l}{2}-1}}\ar@{.}[d]&&\circ\\
\circ&&\circ\ar[d]^{\sigma(v_{1})}\ar@{.}[u]\\
\circ\ar[d]_{u_1}\ar[u]^{v_1}&&\circ\\
\circ&&\circ\ar[u]_{\sigma(u_{1})}\ar@{.}[d]\\
\circ\ar[d]_{u_{\frac{k}{2}}}\ar@{.}[u]&&\circ\\
 \circ\ar[rr]_{b}&&\circ\ar[u]_{\sigma(u_{\frac{k}{2}})}.}
\end{equation}
Let $E^{h}_{i,j}$ be the regular indecomposable representation of dimension $e_{i,j}+h$ which contains $E_{i,j}$ (similarly we define $E^{'h}_{i,j}$). From the regular component of the Auslander-Reiten quiver of $Q'$, we note that the following indecomposable representations $V'\in Rep(Q')$
satisfy property \textit{(Op)} (the other regular indecomposable
representations of $Rep(Q')$ satisfying property \textit{(Op)} are
extensions of these).
\begin{itemize}
\item[(a)] $V_{(0,1)}$ such that $Z'=E^{h}_{1}\oplus E_{0,2}=\nabla Z'$.
\item[(b)] $E_{i-1,j}$, with $1\leq j< i-1\leq l$, such that
$\nabla E_{i-1,j}=E_{i,j+1}$ and
$Z'=E_{i+1,j-1}\oplus E_{i,j}=\nabla Z'$.
\item[(c)] $E_{i-1,j}^h$, with $1\leq i-1< j\leq l$, such that
$\nabla E_{i-1,j}^h=E_{i,j+1}^h$ and
$Z'=E_{i+1,j-1}^h\oplus E_{i,j}^h= \nabla Z'$.
\item[(d)] $E'_{i-1,j}$, with $1\leq i-1<j\leq k+1$, such that
$\nabla E'_{i-1,j}=E'_{i,j+1}$ and
$Z'=E'_{i+1,j-1}\oplus E'_{i,j}=\nabla Z'$.
\item[(e)] $E_{i-1,j}^{'h}$, with $1\leq j\leq i-1\leq k+1$, such that
$\nabla E_{i-1,j}^{'h}=E_{i,j+1}^{'h}$ and
$Z'=E_{i+1,j-1}^{'h}\oplus E_{i,j}^{'h}=\nabla Z'$.
\end{itemize}
If $V$ is the middle term of a short exact sequence $0\rightarrow
V^1\rightarrow V\rightarrow V^2\rightarrow 0$, with $V^1$ and
$V^2$ one of the representations of type (a), (b), (c), (d) or (e), we have
the blocks matrix
$$
Hom_Q(d_{min}^V,\cdot)=\left(\begin{array}{cc}Hom_Q(d_{min}^{V^1},\cdot) & 0 \\
Hom_Q(B,\cdot) & Hom_Q(d_{min}^{V^2},\cdot)\end{array}\right).
$$
where $d_{min}^{V^1}:P^1_1\rightarrow P^1_0$ is the minimal
projective presentation of $V^1$, $d_{min}^{V^2}:P^2_1\rightarrow
P^2_0$ is the minimal projective presentation of $V^2$ and for
some $B\in Hom_Q(P^2_1,P^1_0)$. In general for every blocks matrix
we have $\left(\begin{array}{cc}A & 0\\0 &
C\end{array}\right)=\left(\begin{array}{cc}Id & 0\\-BA^{-1}&
Id\end{array}\right)\cdot \left(\begin{array}{cc}A & 0\\B &
C\end{array}\right)$ if $A$ is invertible. Hence using rows
operations on $Hom_Q(d_{min}^V,\cdot)$, we obtain
$$
Hom_Q(d_{min}^V,\cdot)\approx\left(\begin{array}{cc}Hom_Q(d_{min}^{V^1},\cdot) & 0 \\
0 & Hom_Q(d_{min}^{V^2},\cdot)\end{array}\right).
$$
So it's enough to prove the skew-symmetry of
$Hom_Q(d_{min}^V,\cdot)$ on $SpRep(Q,\beta)$ for $V$ one of
representations of type (a), (b), (c), (d) and (e).\\
Let $\chi$ be the symmetric weight associated to $\alpha$. We
order vertices of $Q$ clockwise from $tb=1$ to $hb=k+l+1$.\\
Let $W\in SpRep(Q,\beta)$. We prove that
$Hom_Q(d_{min}^V,W)$ is skew-symmetric for every regular
indecomposable representation $V$ of type (a), (b), (c), (d) and (e). First
we observe that the associated to $V$ symmetric weight $\chi$ have
components equal to 0, 1 and -1. Let $m_1$ be the first vertex
such that $\chi(m_1)\neq 0$ and $m_s$ the last vertex such that
$\chi(m_s)\neq 0$. Between $-1$ and 1 alternate in correspondence
respectively of sinks and of sources. In particular,
$\chi(m_1)=\pm 1=-\chi(m_s)$ and $\chi(m_i)=1$ or $-1$, for every
$i\in\{2,\ldots,s-1\}$, respectively if $m_i$ is a source or a
sink. We note that, for every $Hom_Q(d_{min}^V,W)$ with $V$ one
representation of type (a), (b), (c), (d) and (e), we can restrict to the
symmetric subquiver of type $A$ which has first vertex $m_1$ and
last vertex $m_s$ and passing through the $\sigma$-fixed vertex of
$Q$. Hence it proceeds as done for type $A$ (see Lemma 3.1 in
\cite{ar}).
\end{proof}
\begin{definizione}
For $V\in Rep(Q,\alpha)$ satisfying property \textit{(Spp)} (resp.
satisfying property \textit{(Op)} such that
$\langle\alpha,\beta\rangle=0$, where $\beta$ is an orthogonal
(resp. symplectic) dimension vector, we define
$$
\begin{array}{rcl}
pf^V:Rep(Q,\beta)&\longrightarrow&\Bbbk\\
W&\longmapsto&c^V(W)=pf(Hom_Q(d^V_{min},W)).
\end{array}
$$
\end{definizione}
\begin{proposizione}\label{SIcx1}
Let $(Q,\sigma)$ be a symmetric quiver of tame type. Let
$\alpha$ be a symmetric dimension vector, $x$ be an admissible
sink and $\varphi^{Sp}_{x,\alpha}$ be as defined in Lemma \ref{kac}.\\
Then $\varphi^{Sp}_{x,\alpha}(c^V)=c^{C^+_{(x,\sigma(x))}V}$ and
$\varphi^{Sp}_{x,\alpha}(pf^W)=pf^{C^+_{(x,\sigma(x))}W}$, where $V$
and $W$ are indecomposables of $Q$ such that $\langle
\underline{dim}\,V,\alpha
\rangle=0=\langle\underline{dim}\,W,\alpha \rangle$ and $W$
satisfies property \textit{(Op)}. In particular
\begin{itemize}
\item[(i)] if $0=\alpha_x\neq\sum_{a\in Q_1:ta=x}\alpha_{ta}$, then
$(\varphi^{Sp}_{x,\alpha})^{-1}(c^{S_x})=0$;
\item[(ii)] if $0\neq\alpha_x=\sum_{a\in Q_1:ta=x}\alpha_{ta}$, then
$\varphi^{Sp}_{x,\alpha}(c^{S_{\sigma(x)}})=0$.
\end{itemize}
\end{proposizione}
\begin{proof} See Proposition 3.3 in \cite{ar}.
\end{proof}
\begin{proposizione}\label{SIcx2}
Let $(Q,\sigma)$ be a symmetric quiver of tame type. Let
$\alpha$ be a symmetric dimension vector, $x$ be an admissible
sink and $\varphi^{O}_{x,\alpha}$ be as defined in Lemma \ref{kac}.\\
Then $\varphi^{O}_{x,\alpha}(c^V)=c^{C^+_{(x,\sigma(x))}V}$ and
$\varphi^{O}_{x,\alpha}(pf^W)=pf^{C^+_{(x,\sigma(x))}W}$, where $V$
and $W$ are indecomposables of $Q$ such that $\langle
\underline{dim}\,V,\alpha
\rangle=0=\langle\underline{dim}\,W,\alpha \rangle$ and $W$
satisfies property \textit{(Spp)}. In particular
\begin{itemize}
\item[(i)] if $0=\alpha_x\neq\sum_{a\in Q_1:ta=x}\alpha_{ta}$, then
$(\varphi^{O}_{x,\alpha})^{-1}(c^{S_x})=0$;
\item[(ii)] if $0\neq\alpha_x=\sum_{a\in Q_1:ta=x}\alpha_{ta}$, then
$\varphi^{O}_{x,\alpha}(c^{S_{\sigma(x)}})=0$.
\end{itemize}
\end{proposizione}
\begin{proof} See Proposition 3.4 in \cite{ar}.
\end{proof}
By Proposition \ref{SIcx1}, Proposition \ref{SIcx2} and by Lemma
\ref{kac}, it follows that if Theorem \ref{tp1} and Theorem
\ref{tp2} are true for $(Q,\sigma)$, then they are true for
$(c_{(x,\sigma(x))}Q,\sigma)$.

\subsection{Semi-invariants of Symmetric quivers of tame type for
dimension vector $ph$} In this section we deal with dimension
vector $ph$ (for definition of $h$, see Proposition \ref{h}). One proves Theorems \ref{tp1} and \ref{tp2} type by
type of quivers $\widetilde{A}$ and $\widetilde{D}$. As said above, we can consider orientations of
symmetric quivers of type $\widetilde{A}$ in Proposition \ref{oA}
and orientation of symmetric quiver $\widetilde{D}^{eq}$. We
proves these theorems only for $\widetilde{A}^{1,1}_{k,l}$. For the
other type the proof is similar (for further details see sec 3.1 in \cite{a1}).\\
We shall call $\Lambda$, $ER\Lambda$ and $EC\Lambda$ respectively
the set of partitions, the set of partition with even rows and the
set of partition of even columns.\\
Theorems \ref{tp1} and \ref{tp2} for type $(1,1,k,l)$ follow from next theorem.
\begin{teorema}
Let $(Q,\sigma)$ be a symmetric quiver of type $(1,1,k,l)$ with
orientation as in (\ref{A11}). Then\\
\textbf{O)} $OSI(Q,ph)$ is generated by the following indecomposable
semi-invariants:\\
if $p$ is even,
\begin{itemize}
\item[a)] $det\,V(u_j)$ with $j\in\{1,\ldots,\frac{k}{2}\}$;
\item[b)] $det\,V(v_j)$ with $j\in\{1,\ldots,\frac{l}{2}\}$;
\item[c)] $pf\,V(b)$
\item[d)] the coefficients $c_i$ of $\varphi^{p-2i}\psi^{2i}$, $0\leq i\leq \frac{p}{2}$, in $det(\psi V(\sigma(\bar{a})\bar{a})+\varphi
V(\bar{b}))$, where $\bar{a}= v_{\frac{l}{2}}\cdots v_1$ and
$\bar{b}=\sigma(u_{1})\cdots \sigma(u_{\frac{k}{2}})b
u_{\frac{k}{2}}\cdots u_1$;
\end{itemize}
if $p$ is odd,
\begin{itemize}
\item[a)] $det\,V(u_j)$ with $j\in\{1,\ldots,\frac{k}{2}\}$;
\item[b)] $det\,V(v_j)$ with $j\in\{1,\ldots,\frac{l}{2}\}$;
\item[c)] the coefficients $c_i$ of $\varphi^{p-2i}\psi^{2i}$, $0\leq i\leq \frac{p-1}{2}$, in $det(\psi V(\sigma(\bar{a})\bar{a})+\varphi
V(\bar{b}))$, where $\bar{a}= v_{\frac{l}{2}}\cdots v_1$ and
$\bar{b}=\sigma(u_{1})\cdots \sigma(u_{\frac{k}{2}})b
u_{\frac{k}{2}}\cdots u_1$.
\end{itemize}
\textbf{Sp)} $SpSI(Q,ph)$ is generated by the following indecomposable
semi-invariants:\\
if $p$ is even,
\begin{itemize}
\item[a)] $det\,V(u_j)$ with $j\in\{1,\ldots,\frac{k}{2}\}$;
\item[b)] $det\,V(v_j)$ with $j\in\{1,\ldots,\frac{l}{2}\}$;
\item[c)] $det\,V(b)$
\item[d)] the coefficients $c_i$ of $\varphi^{p-2i}\psi^{2i}$, $0\leq i\leq \frac{p}{2}$, in $det(\psi V(\sigma(\bar{a})\bar{a})+\varphi
V(\bar{b}))$, where $\bar{a}=v_{\frac{l}{2}}\cdots v_1$ and
$\bar{b}=\sigma(u_{1})\cdots \sigma(u_{\frac{k}{2}})b
u_{\frac{k}{2}}\cdots u_1$;
\end{itemize}
if $p$ is odd, $SpSI(Q,ph)=\Bbbk$.
\end{teorema}
\begin{proof} We proceed by induction on $\frac{k}{2}+\frac{l}{2}$. The
smallest case is $\widetilde{A}^{1,1}_{0,2}$
$$
\xymatrix@-1pc{&2\ar[dr]^{\sigma(a)}&\\
1\ar[rr]_b\ar[ur]^a&&\sigma(1).}
$$
The induction step follows by Lemma \ref{cl} and by Lemma \ref{cls}, so it's enough to prove the theorem for $\widetilde{A}^{1,1}_{0,2}$.\\
\textbf{O)} The ring of orthogonal semi-invariants is
$$
\bigoplus_{{\lambda(a)\in\Lambda \atop\lambda(b)\in EC
\Lambda}}(S_{\lambda(a)}V_1\otimes
S_{\lambda(b)}V_1)^{SL\,V_1}\otimes (S_{\lambda(a)}V_2)^{SO\,V_2}.
$$
By Proposition \ref{i2} we have
\begin{equation}\label{7}
\lambda(a)_j+\lambda(b)_{p-j+1}=k_1
\end{equation}
for every $0\leq j\leq p$ and for some $k_1\in\mathbb{N}$. By Proposition \ref{i3} we have $\lambda(a)=2\mu+(l^p)$ for
some $\mu\in\Lambda$ and for some $l\in\mathbb{N}$. We
consider the summands in which $k_1=1,2$ because the other ones
are generated by products of powers of the generators of this
summands.\\ 
Let $p$ be even. If $k_1=1$ the only solutions of
(\ref{7}) are $\lambda(a)=(1^p)$, $\lambda(b)=0$ and
$\lambda(a)=0$, $\lambda(b)=(1^p)$. Respectively, the summand
$(S_{(1^p)}V_1)^{SL\,V_1}\otimes(S_{(1^p)}V_2)^{SO\,V_2}$ is
generated by a semi-invariant of weight $(1,0)$, i.e
$det\,V(a)=det\,V(\sigma(a))$, and the summand
$(S_{(1^p)}V_1)^{SL\,V_1}$ is generated by a semi-invariant of
weight $(1,0)$, i.e $pf\,V(b)$. If $k_1=2$, the solutions of
(\ref{7}) are $\lambda(a)=(2^{2i})$, $\lambda(b)=(2^{p-2i})$ with
$0\leq i\leq \frac{p}{2}$. So the summand is
$$
\bigoplus_{i=0}^{\frac{p}{2}}(S_{(2^{2i})}V_1\otimes
S_{(2^{p-2i})}V_1)^{SL\,V_1}\otimes (S_{(2^{2i})}V_2)^{SO\,V_2}
$$
which is generated by the coefficients of
$\varphi^{p-2i}\psi^{2i}$ in $det(\psi V(\sigma(a)a)+\varphi
V(b))$, semi-invariants of weight $(2,0)$. In particular for $i=0$
we have
$det\,V(b)$ and for $i=\frac{p}{2}$ we have $det\,V(\sigma(a)a)$.\\
Let $p$ be odd. If $k_1=1$ the only solutions of (\ref{7}) are
$\lambda(a)=(1^p)$, $\lambda(b)=0$. The summand
$(S_{(1^p)}V_1)^{SL\,V_1}\otimes(S_{(1^p)}V_2)^{SO\,V_2}$ is
generated by a semi-invariant of weight $(1,0)$, i.e
$det\,V(a)=det\,V(\sigma(a))$. If $k_1=2$, the solutions of
(\ref{7}) are $\lambda(b)=(2^{2i})$, $\lambda(a)=(2^{p-2i})$ with
$0\leq i\leq \frac{p-1}{2}$. So the summand is
$$
\bigoplus_{i=0}^{\frac{p-1}{2}}(S_{(2^{p-2i})}V_1\otimes
S_{(2^{2i})}V_1)^{SL\,V_1}\otimes (S_{(2^{p-2i})}V_2)^{SO\,V_2}
$$
which is generated by the coefficients of
$\varphi^{2i}\psi^{p-2i}$ in $det(\psi V(\sigma(a)a)+\varphi
V(b))$, semi-invariants of weight $(2,0)$.\\
\textbf{Sp)} The ring of symplectic semi-invariants is
$$
\bigoplus_{{\lambda(a)\in\Lambda \atop\lambda(b)\in ER
\Lambda}}(S_{\lambda(a)}V_1\otimes
S_{\lambda(b)}V_1)^{SL\,V_1}\otimes (S_{\lambda(a)}V_2)^{Sp\,V_2}.
$$
By Proposition \ref{i2} we have
\begin{equation}\label{8}
\lambda(a)_j+\lambda(b)_{p-j+1}=k_1
\end{equation}
for every $0\leq j\leq p$ and for some $k_1\in\mathbb{N}$. By Proposition \ref{i3} we have $\lambda(a)\in EC\Lambda$. We
consider the summands in which $k_1=1,2$ because the other ones
are generated by products of powers of the generators of this
summands.\\ 
Let $p$ be even. If $k_1=1$ the only solutions of
(\ref{8}) are $\lambda(a)=(1^p)$, $\lambda(b)=0$. The summand
$(S_{(1^p)}V_1)^{SL\,V_1}\otimes(S_{(1^p)}V_2)^{Sp\,V_2}$ is
generated by a semi-invariant of weight $(1,0)$, i.e
$det\,V(a)=det\,V(\sigma(a))=pf\,V(\sigma(a)a)$. If $k_1=2$, the
solutions of (\ref{8}) are $\lambda(a)=(2^{2i})$,
$\lambda(b)=(2^{p-2i})$ with $0\leq i\leq \frac{p}{2}$. So the
summand is
$$
\bigoplus_{i=0}^{\frac{p}{2}}(S_{(2^{2i})}V_1\otimes
S_{(2^{p-2i})}V_1)^{SL\,V_1}\otimes (S_{(2^{2i})}V_2)^{Sp\,V_2}
$$
which is generated by the coefficients of
$\varphi^{p-2i}\psi^{2i}$ in $det(\psi V(\sigma(a)a)+\varphi
V(b))$, semi-invariants of weight $(2,0)$. In particular for $i=0$
we have
$det\,V(b)$ and for $i=\frac{p}{2}$ we have $det\,V(\sigma(a)a)$.\\
If $p$ is odd there not exist any non-trivial symplectic
representations because a symplectic space of dimension odd
doesn't exist. So we have $SpSI(Q,ph)=\Bbbk$.\end{proof}

From previous theorem, theorems \ref{tp1} and \ref{tp2} follow for symmetric quivers of type $(1,1,k,l)$.
\begin{proof}[Proof of Theorems \ref{tp1} and \ref{tp2} for type $(1,1,k,l)$]
First of all we note that, by definition of $c^W$ and $pf^W$, when
we have it, are not zero if $0=\langle
\underline{dim}\,W,ph\rangle=p\langle
\underline{dim}\,W,h\rangle=-p\langle
h,\underline{dim}\,W\rangle$, so we have to consider only regular
representations $W$. Moreover it is enough to consider only simple
regular representations $W$, because the other regular
representations are extensions of simple regular ones and so, by
Lemma \ref{cVcV'}, we obtain the $c^W$ and $pf^W$ with non-simple
regular $W$ as products of those with  simple regular $W$. Now we
check only for $\widetilde{A}^{1,1}_{k,l}$ that the generators
found, in previous theorem, for $SpSI(Q,ph)$ and $OSI(Q,ph)$ are of type $c^{W}$, for
some simple regular $W$, and $pf^W$, for some simple regular $W$
satisfying property \textit{(Op)} in symplectic case and
\textit{(Spp)} in orthogonal case (see Lemma \ref{ss1}).
\begin{itemize}
\item[\textbf{Sp)}] Let $V$ be a symplectic representation. We recall that, in this case, $p$ has to be even. The minimal projective resolution of $E_0$ is
$$
0\longrightarrow
P_{h(\sigma(v_1))}\stackrel{d^{E_0}_{min}}{\longrightarrow}P_{t(\sigma(v_1))}\longrightarrow
E_0\longrightarrow 0,
$$
where $d^{E_0}_{min}=\sigma(v_1)$. So we have
$c^{E_0}(V)=det(V(\sigma(v_1)))=det(V(v_1))=c^{E_1}(V)$.
Similarly, we obtain
$c^{E_i}(V)=det(V(v_i))=det(V(\sigma(v_i)))=c^{E_{l-i+1}}(V)$
for every $i\in\{2,\ldots,l\}$,
$c^{E'_0}(V)=det(V(\sigma(u_1)))=det(V(u_1))=c^{E'_1}(V)$,
$c^{E'_i}(V)=det(V(u_i))=det(V(\sigma(u_i)))=c^{E'_{k-i+2}}(V)$
for every $i\in\{2,\ldots,k\}\setminus\{\frac{k}{2}+1\}$,
$c^{E'_{\frac{k}{2}+1}}(V)=det(V(b))$ and
$c^{V_{(\varphi,\psi)}}(V)=det(\psi
V(\sigma(\bar{a})\bar{a})+\varphi V(\bar{b}))$;
\item[\textbf{O)}] if $V$ is an orthogonal representation, the only difference with the symplectic case is,
when $p$ is even, we have $pf^{E'_{\frac{k}{2}+1}}(V)=pf(V(b))$
and $pf^{V_{(\varphi,\psi)}}(V)=pf(\psi V(a)+\varphi V(b))$, in
fact $E'_{\frac{k}{2}+1}$ satisfies property \textit{(Spp)}, since $C^+E'_{\frac{k}{2}+1}=E'_{\frac{k}{2}+2}=E'_{\sigma(\frac{k}{2}+1)}=\nabla E'_{\frac{k}{2}+1}$ and the almost split sequence
$$
0\rightarrow E'_{\frac{k}{2}+2}\rightarrow Z\rightarrow E'_{\frac{k}{2}+1}\rightarrow 0 
$$
has $Z= E'_{\frac{k}{2}+2,\frac{k}{2}+1}\in SpRep(Q)$.
\end{itemize}
\end{proof}

\subsection{Semi-invariants of symmetric quivers of tame type for any regular dimension vector}
In this
section we prove theorems \ref{tp1} and \ref{tp2} for symmetric
quivers of tame type and any regular symmetric dimension vector
$d$.\\
We will use the same notation of section 3.1. For the type
$\widetilde{A}$ we call $a_0=tv_1=tu_1$, $x_i=hv_i$ for every
$i\in\{1,\dots,\frac{l}{2}\}$ and $y_i=hv_i$ for every
$i\in\{1,\dots,\frac{k}{2}\}$. For the type $\widetilde{D}$ we
call $t_1=ta$, $t_2=tb$ and $z_i=tc_i$ for every $i$ such that
$c_i\in (Q_1^+\sqcup Q_1^{\sigma})\setminus\{a,b\}$.\\
First we consider the decomposition of $d$ for the
symmetric quivers.\\
Let $(Q,\sigma)$ be a symmetric quiver of tame type and let
$\Delta=\{e_i|\,i\in I=\{0,\ldots,u\}\}$, $\Delta'=\{e'_i|\,i\in
I'= \{0,\ldots,v\}\}$ and $\Delta''=\{e''_i|\,i\in
I''=\{0,\ldots,w\}\}$ be the three $C^+$-orbits of
nonhomogeneous simple regular representations of the underlying
quiver $Q$ (see Proposition \ref{orbitetau}).\\
We shall call $I_{\delta}=\{i\in I|\,e_i=\delta e_i\}$
(respectively $I'_{\delta}$ and $I''_{\delta}$).
Let $[x]:=max\{z\in \mathbb{N}|z\leq x\}$ be the floor of $x\in\mathbb{R}$.
\begin{lemma}\label{dI}
\begin{itemize}
\item[(a)] For symmetric quivers of tame type with central symmetry, we have $\Delta=\Delta'$ and so $I=I'$.
\item[(b)] For symmetric quivers of tame type with symmetry respect to a central vertical line, we have decomposition $I=I_+\sqcup I_{\delta}\sqcup I_-$ such that 
\begin{itemize}
\item[(i)] $I_+=\{2,\ldots,[\frac{u}{2}]+2\}$, $I_{\delta}=\empty$, $I_-=I\setminus I_+$;
\item[(ii)] $I_+=\{2,\ldots,[\frac{u}{2}]+1\}$, $I_{\delta}=\{1\}$, $I_-=I\setminus (I_+\sqcup I_{\delta})$;
\item[(iii)] $I_+=\{2,\ldots,[\frac{u}{2}]+1\}$, $I_{\delta}=\{1,[\frac{u}{2}]+2\}$, $I_-=I\setminus (I_+\sqcup I_{\delta})$.
\end{itemize}
\end{itemize}
Similarly for $I'$ but replacing $v$ with $u$. $I''=I_+''\sqcup I_-''$ such that $I_+''=\{2\}$ and $I_-''=I''\setminus I_+''$. 
\end{lemma}
\begin{proof}
It follows by section 6 pages 40 and 46 in \cite{dr} and by definition of $\delta$ (Definition \ref{delta}).
\end{proof}
\begin{proposizione}\label{dc}
Let $(Q,\sigma)$ be a symmetric quiver of tame type and let $I_+$,
$I_{\delta}$, $I'_+$, $I'_{\delta}$ and $I''_+$ be
as above. Any regular symmetric dimension vector can be written
uniquely in the following form:
\begin{equation}\label{fdc}
d=ph+\sum_{i\in I_+}p_i(e_i+\delta e_i)+\sum_{i\in I_{\delta}}p_i
e_i+\sum_{i\in I'_+}p'_i(e'_i+\delta e'_i)+\sum_{i\in
I'_{\delta}}p'_i e'_i+\sum_{i\in I''_+}p''_i(e''_i+\delta e''_i)
\end{equation}
for some non-negative $p,p_i,p'_i,p''_i$ with at least one
coefficient in each family $\{p_i|\;i\in I_+\sqcup I_{\delta}\}$,
$\{p'_i|\;i\in I_+'\sqcup I_{\delta}'\}$, $\{p''_i|\;i\in I''_+\}$
being zero. In particular, in the symplectic case,
\begin{itemize}
\item[i)] if $Q$ has one $\sigma$-fixed vertex and one
$\sigma$-fixed arrow (i.e. $Q=\widetilde{A}_{k,l}^{1,1}$), then
$p_{\frac{l}{2}+1}$ and $p'_1$ have to be even,
\item[ii)] if $Q$ has one or two $\sigma$-fixed vertices and it
has not any $\sigma$-fixed arrows (i.e.
$Q=\widetilde{A}_{k,l}^{0,2}$ or $\widetilde{D}_{n}^{0,1}$), then
both $p_i$'s and $p'_j$'s, with $i\in I_{\delta}$ and $j\in
I'_{\delta}$, have to be even.
\end{itemize}
\end{proposizione}
\begin{proof} It follows by Lemma \ref{dI} and by decomposition
of any regular dimension vector of the underlying quiver of
$(Q,\sigma)$. In particular, since symplectic spaces with odd
dimension don't exist, it implies \textit{i)} and
\textit{ii)}. \end{proof}

Graphically we can represent, in the symmetric case, $\Delta$ (similarly $\Delta'$ and
$\Delta''$) as the polygons
$$
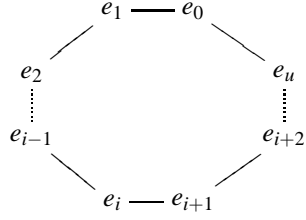
\xymatrix@-1pc{
&e_1\ar@{-}[r]\ar@{-}[dl]&e_0\ar@{-}[dr]&\\
e_2\ar@{.}[d]&&&e_u\ar@{.}[d]\\
e_{i-1}\ar@{-}[dr]&&&\ar@{-}[dl]e_{i+2}\\
&e_{i}\ar@{-}[r]&e_{i+1}&}
$$
if $Q=\widetilde{A}_{k,k}^{0,0}$ and
$$
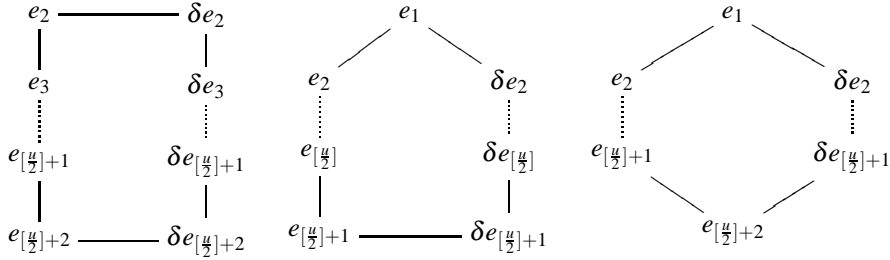
\begin{array}{ccc}
\xymatrix@-1pc{
 e_2\ar@{-}[rr]\ar@{-}[d]&&\delta e_2\ar@{-}[d]\\
 e_3\ar@{.}[d]&&\delta e_3\ar@{.}[d]\\
 e_{[\frac{u}{2}]+1}\ar@{-}[d]&&\delta e_{[\frac{u}{2}]+1}\ar@{-}[d]\\
  e_{[\frac{u}{2}]+2}\ar@{-}[rr]&&\delta e_{[\frac{u}{2}]+2}}&
\xymatrix@-1pc{
&e_1\ar@{-}[dr]\ar@{-}[dl]&\\
 e_2\ar@{.}[d]&&\delta e_2\ar@{.}[d]\\
 e_{[\frac{u}{2}]}\ar@{-}[d]&&\delta e_{[\frac{u}{2}]}\ar@{-}[d]\\
  e_{[\frac{u}{2}]+1}\ar@{-}[rr]&&\delta e_{[\frac{u}{2}]+1}}&
\xymatrix@-1pc{
&e_1\ar@{-}[dr]\ar@{-}[dl]&\\
 e_2\ar@{.}[d]&&\delta e_2\ar@{.}[d]\\
e_{[\frac{u}{2}]+1}\ar@{-}[dr]&&\delta e_{[\frac{u}{2}]+1}\ar@{-}[dl]\\
&e_{[\frac{u}{2}]+2}&}
 \end{array}
$$
with a reflection respect to a central vertical line, in the other
cases.
\begin{definizione}
We define an involution $\sigma_I$ on the set of indices $I$ such
that $e_{\sigma_I(i)}=\delta e_i$ for every $i\in I$. Hence
$\sigma_I(I)=I'$ for $\widetilde{A}_{k,k}^{0,0}$ and
$\sigma_{I}I_+=I_-$, $\sigma_I I_{\delta}=I_{\delta}$ for the
other cases. Similarly we define an involution $\sigma_{I'}$ and
an involution $\sigma_{I''}$ respectively on $I'$ and on $I''$.
\end{definizione}
In the remainder of the section, we shall call
\begin{equation}\label{d'}
d'=\sum_{i\in I_+}p_i(e_i+\delta e_i)+\sum_{i\in
I_{\delta}}p_i e_i+\sum_{i\in I'_+}p'_i(e'_i+\delta
e'_i)+\sum_{i\in I'_{\delta}}p'_i e'_i+\sum_{i\in
I''_+}p''_i(e''_i+\delta e''_i).
\end{equation}
\begin{proposizione}\label{oa}
If $d$ is regular with decomposition (\ref{fdc}) such that $d= d'$
or $d$ is not regular then $SpRep(Q, d)$ (respectively
$ORep(Q,d)$) has an open $Sp(Q, d)$-orbit (respectively
$O(Q,d)$-orbit).
\end{proposizione}
\begin{proof} If $d=d'$, we have no indecomposable of dimension
vector $ph$ and so there are finitely many orbits.
If $d$ is not regular, it follows from \cite{r2}, Theorem 3.2. \end{proof}
In the next, $d$ shall be a regular symmetric dimension vector with
decomposition (\ref{fdc}) with $p\geq 1$ and $p\neq 0$. Now we
shall describe the generators of $SpSI(Q,d)$ and $OSI(Q,d)$.\\
By Proposition \ref{oa}, $Sp(Q,d')$ acting
on $SpRep(Q,d')$ has an open orbit
so, by Lemma \ref{sk}, dimension of $SpSI(Q,d')_{\chi'}$ is 0 or 1. This
allows us to identify one non-zero element of $SpSI(Q,d)_{\chi}$ with the element of
$SpSI(Q,ph)_{\chi}$ to which
it restricts. Similarly one proceeds for $OSI(Q,d)$.\\
We proceed now to describe the generators of the algebra $SpSI(Q,
d)$ (respectively $OSI(Q,d)$). If the corresponding $I,\,I',\,I''$
are not empty, we label the vertices $e_i,\, e'_i,\, e''_i$ of the
polygons $\Delta,\, \Delta',\, \Delta''$ with the coefficients
$p_i,\, p'_i,\, p''_i$. We recall that
\begin{itemize}
\item[a)] we have to label with $p_i$ (respectively with $p'_i$ and $p''_i$)
both vertices $e_i$ and $\delta e_i$, i.e $p_i=p_{\sigma_I(i)}$
(respectively $p'_i=p'_{\sigma_I'(i)}$ and
$p''_i=p''_{\sigma_I''(i)}$), if $e_i\neq\delta e_i$.
\end{itemize}
and in the symplectic case, by \textit{i)} and \textit{ii)} of
Proposition \ref{dc}
\begin{itemize}
\item[b)] for $\widetilde{A}_{k,l}^{1,1}$, $p_{\frac{l}{2}+1}$ and $p'_1$
have to be even,
\item[c)] for $\widetilde{A}_{k,l}^{0,2}$ and
$\widetilde{D}_{n}^{0,1}$,  $p_i\in I_{\delta}$ and $p'_i\in
I'_{\delta}$ have to be even.
\end{itemize}
We shall call these labelled polygons respectively $\Delta(d),\,
\Delta'(d),\, \Delta''(d)$.
\begin{definizione}
We shall say that the labelled arc
$\xymatrix@-1pc{p_i\ar@{-}[r]&\ar@{.}[r]&\ar@{-}[r]&p_j}$ (in
clockwise orientation) of the labelled polygon $\Delta(d)$ is
admissible if $p_i=p_j$ and $p_i<p_k$ for every its interior
labels $p_k$. We denote such a labelled arc
$\xymatrix@-1pc{p_i\ar@{-}[r]&\ar@{.}[r]&\ar@{-}[r]&p_j}$ by
$[i,j]$, and we define $p_i=p_j$ the index $ind[i,j]$ of $[i,j]$.
Similarly we define admissible arcs and their indexes for the
labelled polygons $\Delta'(d)$ and $\Delta''(d)$.
\end{definizione}
We denote by $\mathcal{A}(d)$, $\mathcal{A}'(d)$,
$\mathcal{A}''(d)$ the sets of all admissible labelled arcs in the
polygons $\Delta(d)$, $\Delta'(d)$, $\Delta''(d)$ respectively. In
particular we note that if $d=ph$, then the polygons $\Delta(d)$,
$\Delta'(d)$, $\Delta''(d)$ are labelled by zeros and so
$\mathcal{A}(d)$, $\mathcal{A}'(d)$, $\mathcal{A}''(d)$ consist of
all edges of respective polygons. With these notations we have the
following
\begin{proposizione}\label{fij}
For each arc $[i,j]$ from $\mathcal{A}(d)$ (respectively
$\mathcal{A}'(d)$ and $\mathcal{A}''(d)$) there exists in
$SpSI(Q,d)$ and in $OSI(Q,d)$ a non zero semi-invariant
\begin{itemize}
\item[(i)] of type $c^{E_{i-1,j}}$ (respectively $c^{E'_{i-1,j}}$
and $c^{E''_{i-1,j}}$) or of type $c^{V_{(\varphi,\psi)}}$, with
$(\varphi,\psi)\in\{(1,0),(0,1),(1,1)\}$;
\item[(ii)] of type $pf^{E_{i-1,j}}$ (respectively $pf^{E'_{i-1,j}}$
and $pf^{E''_{i-1,j}}$) or of type $pf^{V_{(\varphi,\psi)}}$, with
$(\varphi,\psi)\in\{(1,0),(0,1),(1,1)\}$, if $E_{i,j-1}$,
$E'_{i-1,j}$, $E''_{i-1,j}$ and $V_{(\varphi,\psi)}$ satisfy
property \textit{(Op)} in the symplectic case and property
\textit{(Spp)} in the orthogonal case.
\end{itemize}
\end{proposizione}
Let $c_0,\ldots,c_t$ be the coefficients of $\varphi^{t-i}\psi^i$ in $c^{V_{(\varphi,\psi)}}$ or $pf^{V_{(\varphi,\psi)}}$. We note case by case, from previous section, that  $t$ can be $\frac{p-1}{2}$, $\frac{p}{2}$ or
$p$. The generators of
algebras $SpSI(Q,d)$ and $OSI(Q,d)$ are described by the following
theorem
\begin{teorema}\label{stq}
Let $(Q,\sigma)$ a symmetric quiver of tame type and $d=ph+d'$ the
decomposition of a regular symmetric dimension vector $d$ with
$p\geq 1$. Then $SpSI(Q,d)$ (respectively $OSI(Q,d)$) is generated
by
\begin{itemize}
\item[(i)] $c_0,\ldots,c_t$;
\item[(ii)] $c^{E_{i-1,j}}$, $c^{E'_{r-1,s}}$, $c^{E''_{f-1,g}}$ and $c^{V_{(\varphi,\psi)}}$ with $[i,j]\in\mathcal{A}(d)$,
$[r,s]\in\mathcal{A}'(d)$, $[f,g]\in\mathcal{A}''(d)$ and
$(\varphi,\psi)\in\{(1,0),(0,1),(1,1)\}$;
\item[(iii)] $pf^{E_{i-1,j}}$, $pf^{E'_{r-1,s}}$, $pf^{E''_{f-1,g}}$ and $pf^{V_{(\varphi,\psi)}}$ with $[i,j]\in\mathcal{A}(d)$,
$[r,s]\in\mathcal{A}'(d)$, $[f,g]\in\mathcal{A}''(d)$ and
$(\varphi,\psi)\in\{(1,0),(0,1),(1,1)\}$, if $E_{i-1,j}$,
$E'_{r-1,s}$, $E''_{f-1,g}$ and $V_{(\varphi,\psi)}$ satisfy
property \textit{(Op)} (respectively property \textit{(Spp)}).
\end{itemize}
\end{teorema}
For every regular dimension vector $d$, we note that $\langle h,d\rangle=0$ and
$$
\langle\underline{dim} E_{i-1,j},d\rangle=0\Leftrightarrow
p_i=p_j.
$$
So theorem \ref{stq} is equivalent to theorems \ref{tp1} and
\ref{tp2}.
\subsubsection{Generic decomposition for symmetric quivers}
In the
proof of theorem \ref{stq}, we use the notion of generic decomposition of the symmetric
dimension vector $d$ (see remark 2.8(a) in \cite{k1}).
\begin{definizione}
A decomposition $\alpha=\beta_1\oplus\cdots\oplus\beta_q$ of a
dimension vector $\alpha$ is called generic if there is a Zariski
open subset $\mathcal{U}$ of $Rep(Q,\alpha)$ such that each
$U\in\mathcal{U}$ decomposes in $U=\bigoplus_{i=1}^q U_i$ with
$U_i$ indecomposable representation of dimension $\beta_i$, for
every $i\in\{1,\ldots,q\}$.
\end{definizione}
\begin{definizione}
\begin{itemize}
\item[(1)] A decomposition $\alpha=\beta_1\oplus \cdots\oplus\beta_q$ of a symmetric
dimension vector $\alpha$ is called symplectic generic if there is
a Zariski open subset $\mathcal{U}$ of $SpRep(Q,\alpha)$ such that
each $U\in\mathcal{U}$ decomposes in $U=\bigoplus_{i=1}^q U_i$
with $U_i$ indecomposable symplectic representation of dimension
$\beta_i$, for every $i\in\{1,\ldots,q\}$.
\item[(2)] A decomposition $\alpha=\beta_1\oplus\cdots\oplus\beta_q$ of a symmetric
dimension vector $\alpha$ is called orthogonal generic if there is
a Zariski open subset $\mathcal{U}$ of $ORep(Q,\alpha)$ such that
each $U\in\mathcal{U}$ decomposes in $U=\bigoplus_{i=1}^q U_i$
with $U_i$ indecomposable orthogonal representation of dimension
$\beta_i$, for every $i\in\{1,\ldots,q\}$.
\end{itemize}
\end{definizione}
For tame quivers the generic decomposition of any regular
dimension vector is given by
results of section 3 in \cite{r2}.\\
We describe this decomposition explicitly for a symmetric regular
dimension vector $d$ with decomposition (\ref{fdc}).\\
In the remainder of this section we set
\begin{equation}\label{d'1}
\bar{d}=\sum_{i\in I_+}p_i(e_i+\delta e_i)+\sum_{i\in
I_{\delta}}p_i e_i,
\end{equation}
\begin{equation}
\bar{d}'=\sum_{i\in I'_+}p'_i(e'_i+\delta e'_i)+\sum_{i\in
I'_{\delta}}p'_i e'_i
\end{equation}
\begin{equation}
\bar{d}''=\sum_{i\in I''_+}p''_i(e''_i+\delta e''_i).
\end{equation}
\begin{oss}
We can assume $p_i=0$ for $i\in I_{\delta}$ or $p_{i}=0$, for $i\in I_+$, and so
$p_{\sigma_I(i)}=0$.
\end{oss}
\begin{definizione}
We divide the polygon $\Delta(\bar{d})$ in two parts:
\begin{itemize}
\item[(i)] the up part $\Delta_{up}(\bar{d})$ is the part of $\Delta(\bar{d})$ from $p_{i-1}$ to
$p_{\sigma_I(i-1)}$;
\item[(ii)] the down part $\Delta_{down}(\bar{d})$ is the part of $\Delta(\bar{d})$ from $p_{\sigma_I(i+1)}$ to $p_{i+1}$.
\end{itemize}
Similarly for $\Delta'$ and $\Delta''$.
\end{definizione}
\begin{oss}
We note that if $p_i=0$ with $i\in I_{\delta}$, then we have only
the part $\Delta_{up}$ or the part $\Delta_{down}$.
\end{oss}
We consider $\Delta$, similarly one proceeds for $\Delta'$ and
$\Delta''$.
\begin{definizione}
We shall call symmetric arc, an arc invariant under $\sigma_I$,
i.e. an arc of type $[i,\sigma_I(i)]$.
\end{definizione}
\begin{oss}
By the division of $\Delta$ in $\Delta_{up}$ and $\Delta_{down}$,
we note that all symmetric arcs pass through the same
$\sigma_I$-fixed vertex of $\Delta$ or through the same
$\sigma_I$-fixed edge of $\Delta$.
\end{oss}
\begin{lemma}\label{3218}
Let $(Q,\sigma)$ be a symmetric quiver of tame type.
\begin{itemize}
\item[(i)] If $n=\sigma_I(n)$ then either there exists unique $x\in
Q_0^{\sigma}$ such that $e_n(x)\neq 0$ or there exists unique
$a\in Q_1^{\sigma}$ such that $e_n(ta)\neq 0$.
\item[(ii)] If $\xymatrix@-1pc{n\ar@{-}[r]&\sigma_I(n)}$ is
a $\sigma_I$-fixed edge in $\Delta$, then there exists unique
$a\in Q_1^{\sigma}$ such that $e_n(ta)\neq  0$.
\end{itemize}
\end{lemma}
\begin{proof}
One proves type by type, using Lemma \ref{dI} (for further details, see Lemma 3.2.18 in \cite{a1}).
\end{proof}
\begin{definizione}
\begin{itemize}
\item[(i)] If $n=\sigma_I(n)$, we call $x(n)$ the unique $x\in
Q_0^{\sigma}$ such that $e_n(x)\neq 0$.
\item[(ii)] If $n=\sigma_I(n)$ or $\xymatrix@-1pc{n\ar@{-}[r]&\sigma_I(n)}$ is
a $\sigma_I$-fixed edge in $\Delta$, we call $a(n)$ the unique
$a\in Q_1^{\sigma}$ such that $e_n(ta)\neq  0$.
\end{itemize}
\end{definizione}
\begin{definizione}
For every arc $[i,j]$ in $\Delta$, we define
$$
e_{[i,j]}=\sum_{k\in [i,j]}e_k.
$$
\end{definizione}
\begin{definizione}
\begin{itemize}
\item[(i)] $\mathcal{A}_+(\bar{d}):=\{[i,j]\in\mathcal{A}(\bar{d})|\;[i,j]\subset I_+\}$
\item[(ii)] $\mathcal{A}^k_+(\bar{d}):=\{[i,j]\in\mathcal{A}(\bar{d})|\;[i,j]\subset I_+,ind[i,j]=k\}$.
\item[(iii)] $\mathcal{A}^k_{\sigma_I}(\bar{d})=\{[i,j]=\sigma_I[i,j]\in\mathcal{A}(\bar{d})|\;ind[i,j]=k\}$.
\end{itemize}
\end{definizione}
\begin{oss}
$[i,j]\subset I_+$ if and only if
$[\sigma_I(j),\sigma_I(i)]\subset I_{-}$ and
$ind[i,j]=ind[\sigma_I(j),\sigma_I(i)]$.
\end{oss}
First we consider all the admissible arcs in
$\mathcal{A}^r_{\sigma_I}(\bar{d})\cup \mathcal{A}^r_{+}(\bar{d})$
such that $r=max\{p_k\}$. So we get
$$
\sum_{i\in I_+}p_i(e_i+\delta e_i)+\sum_{i\in
I_{\delta}}p_ie_i=
$$
\begin{equation}\label{barp}
\sum_{i\in I_+}\tilde{p}_i(\tilde{e}_i+\delta
\tilde{e}_i)+\sum_{i\in I_{\delta}}\tilde{p}_i\tilde{e}_i
+\left(\bigoplus_{[i,j]\in
\mathcal{A}^r_{+}(\bar{d})}(e_{[i,j]}+\delta
e_{[i,j]})+\bigoplus_{[i,\sigma_I(i)]\in
\mathcal{A}^r_{\sigma_I}(\bar{d})}e_{[i,\sigma_I(i)]}\right),
\end{equation}
where $max(\tilde{p}_i)=r-1$. Then we repeat the procedure for
(\ref{barp}) and so on we have
$$
\sum_{i\in I_+}p_i(e_i+\delta e_i)+\sum_{i\in I_{\delta}}p_ie_i=
$$
\begin{equation}\label{vj,k}
\bigoplus_{k=1}^{r}\left(\bigoplus_{[i,j]\in
\mathcal{A}^k_{+}(\bar{d})}(e_{[i,j]}+\delta
e_{[i,j]})+\bigoplus_{[i,\sigma_I(i)]\in
\mathcal{A}^k_{\sigma_I}(\bar{d})}e_{[i,\sigma_I(i)]}\right).
\end{equation}
\begin{oss}
\begin{itemize}
\item[(i)] If $[i,j]$ and $[i',j']$ are two admissible arcs in $\mathcal{A}(\bar{d})$ such that $[i,j]\supseteq [i',j']$, then
$ind[i,j]\leq ind[i',j']$.
\item[(ii)] If there not exists $[i,j]\in\mathcal{A}^h_{\sigma_I}(\bar{d})\cup
\mathcal{A}^h_{+}(\bar{d})$ such that $[i,j]\supseteq [i',j']$ for
some $[i',j']\in\mathcal{A}^k_{\sigma_I}(\bar{d})\cup
\mathcal{A}^k_{+}(\bar{d})$, then the symmetric dimension vector
corresponding to $[i',j']$ appears $k$-times in the decomposition
(\ref{vj,k}), with $1\leq h< k$.
\end{itemize}
\end{oss}
\begin{definizione}
Let $[i_1,j_1],\ldots,[i_k,j_k]$ be the admissible arcs such that
$[i_1,j_1]\supseteq\cdots\supseteq[i_k,j_k]$, with $k\geq 1$. We
define $q_{[i_h,j_h]}=ind[i_h,j_h]-ind[i_{h-1},j_{h-1}]$ for every
$1\leq h\leq k$, where $ind[i_0,j_0]=0$.
\end{definizione}
We note that for every $[i,j]\in
\mathcal{A}^k_{\sigma_I}(\bar{d})\cup \mathcal{A}^k_{+}(\bar{d})$,
$q_{[i,j]}$ is the multiplicity of the symmetric dimension vector
corresponding to $[i,j]$ in the
decomposition (\ref{vj,k}).\\
Finally we have
$$
\sum_{i\in I_+}p_i(e_i+\delta e_i)+\sum_{i\in I_{\delta}}p_ie_i=
$$
\begin{equation}\label{qi,j}
\bigoplus_{[i,j]\in\mathcal{A}_+(\bar{d}) }(e_{[i,j]}+\delta
e_{[i,j]})^{\oplus q_{[i,j]}} +
\bigoplus_{[i,\sigma_I(i)]\in\mathcal{A}(\bar{d})
}(e_{[i,\sigma_I(i)]})^{\oplus q_{[i,\sigma_I(i)]}}.
\end{equation}
\begin{esempi}\label{esempioe}
If $\Delta$ is of the form
\begin{equation}\label{esempio}
\xymatrix@-1pc{&e_1=\delta e_1\ar@{-}[dl]\ar@{-}[dr]&\\
 e_2\ar@{-}[d]&&\delta e_2=e_{\sigma_I(2)}\ar@{-}[d]\\
 e_3\ar@{-}[dr]&&\delta e_3=e_{\sigma_I(3)}\ar@{-}[dl]\\
&e_4=\delta e_4&}
\end{equation}
and $p_1=4$, $p_2=3$, $p_3=0$ and $p_4=2$, then
$[2,\sigma_I(2)]=\{2,1,\sigma_I(2)\}\subset I_+\sqcup
I_{\delta}\sqcup I_-$ with
$q_{[2,\sigma_I(2)]}=ind[2,\sigma_I(2)]=3$, $[1,1]=\{1\}\in I_+$
with $q_{[1,1]}=ind[1,1]-ind[2,\sigma_I(2)]=1$ and $[4,4]=\{4\}\in
I_{\delta}$ with $q_{[4,4]}=ind[4,4]=2$. So we have
$$
\sum_{i\in I_+}p_i(e_i+\delta e_i)+\sum_{i\in I_{\delta}}p_ie_i=
((e_2+\delta e_2)+e_1)^{\oplus 3}\oplus e_1\oplus
(e_4)^{\oplus 2}.
$$
\end{esempi}
Similarly we proceed with the decomposition of $\bar{d}'$ and
$\bar{d}''$. So we have the following
\begin{proposizione}
Let $(Q,\sigma)$ be a symmetric quiver of tame type and let $d$ be
a symmetric dimension vector of a representation of the underlying
quiver $Q$ with decomposition (\ref{fdc}). Then
$$
d=\bigoplus_{i=1}^p h+\bigoplus_{[i,j]\in\mathcal{A}_+(\bar{d})
}(e_{[i,j]}+\delta e_{[i,j]})^{\oplus q_{[i,j]}} +
\bigoplus_{[i,\sigma_I(i)]\in\mathcal{A}(\bar{d})}(e_{[i,\sigma_I(i)]})^{\oplus
q_{[i,\sigma_I(i)]}}+
$$
$$
\bigoplus_{[i,j]\in\mathcal{A}'_+(\bar{d}')}(e'_{[i,j]}+\delta
e'_{[i,j]})^{\oplus q'_{[i,j]}}
+\bigoplus_{[i,\sigma_{I'}(i)]\in\mathcal{A}'(\bar{d}')
}(e'_{[i,\sigma_{I'}(i)]})^{\oplus q'_{[i,\sigma_{I'}(i)]}}+
$$
\begin{equation}\label{dgvs}
\bigoplus_{[i,j]\in\mathcal{A}''_+(\bar{d}'')}(e''_{[i,j]}+\delta
e''_{[i,j]})^{\oplus q''_{[i,j]}}
+\bigoplus_{[i,\sigma_{I''}(i)]\in\mathcal{A}''(\bar{d}'')
}(e''_{[i,\sigma_{I''}(i)]})^{\oplus q''_{[i,\sigma_{I''}(i)]}}
\end{equation}
is the generic decomposition of $d$.
\end{proposizione}
We restrict to regular symplectic and orthogonal dimension vectors. We
modify generic decomposition (\ref{dgvs}) of $d=(d_i)_{i\in Q_0}$
to get symplectic generic decomposition of $d$ or orthogonal
generic decomposition of $d$.\\
Let $[i,j]$ be an arc in $\Delta_{up}$ and let $[h,k]$ be an arc
in $\Delta_{down}$. If $E_{[i,j]}$ is the regular indecomposable
symplectic (respectively orthogonal) representation of
$(Q,\sigma)$ corresponding to $[i,j]$ and $E_{[h,k]}$ is the
regular indecomposable symplectic (respectively orthogonal)
representation of $(Q,\sigma)$ corresponding to $[h,k]$, then
$$
Hom_Q(E_{[i,j]},E_{[h,k]})=0=Hom_Q(E_{[h,k]},E_{[i,j]})
$$
and
$$
Ext_Q^1(E_{[i,j]},E_{[h,k]})=0=Ext_Q^1(E_{[h,k]},E_{[i,j]}).
$$
So we deal separately with $\Delta_{up}$ and $\Delta_{down}$. We
consider $I=I^{up}\sqcup I^{down}$, $I_+=I_+^{up}\sqcup
I_+^{down}$ and $I_{\delta}=I_{\delta}^{up}\sqcup
I_{\delta}^{down}$. We have the decomposition
$\bar{d}=\bar{d}_{up}+\bar{d}_{down}$, where
\begin{equation}\label{d'1up}
\bar{d}_{up}=\sum_{i\in I^{up}_+}p_i(e_i+\delta e_i)+\sum_{i\in
I^{up}_{\delta}}p_i e_i
\end{equation}
and
\begin{equation}
\bar{d}_{down}=\sum_{i\in I^{down}_+}p_i(e_i+\delta
e_i)+\sum_{i\in I^{down}_{\delta}}p_i e_i.
\end{equation}
By what has be said, the symplectic (respectively orthogonal)
generic decomposition of $\bar{d}$ is direct sum of the symplectic
(respectively orthogonal) generic decomposition of $\bar{d}_{up}$
and the symplectic (respectively orthogonal) generic decomposition
of $\bar{d}_{down}$.
\begin{oss}
\begin{itemize}
\item[(i)] In the symplectic case, since $\bar{d}_{x}$ has to
be even for every $x\in Q_0^{\sigma}$, we have to modify the
symmetric dimension vectors corresponding to the arcs passing
through the $\sigma_I$-fixed vertex $n$ such that there exists
$x=x(n)\in Q_0^{\sigma}$.
\item[(ii)] In the orthogonal case, we have to modify the symmetric
dimension vectors corresponding to the arcs passing through the
$\sigma_I$-fixed vertex $n$ such that $\bar{d}_{ta(n)}$ is even
and those corresponding to the arcs passing through the
$\sigma_I$-fixed edge $\xymatrix@-1pc{n\ar@{-}[r]&\sigma_I(n)}$
such that $\bar{d}_{ta(n)}$ is even.
\item[(iii)] We have to modify also $ph+e_{[i,\sigma_I(i)]}$, with
$p$ odd, if $[i,\sigma_I(i)]$ is like in part (i) (respectively
part (ii)), since $h+e_{[i,\sigma_I(i)]}$ is the dimension vector
of regular indecomposable symplectic (respectively orthogonal)
representation.
\end{itemize}
\end{oss}
\begin{definizione}
\begin{itemize}
\item[(i)] $\mathcal{A}^{up}(\bar{d})=\{[i,j]\in
\mathcal{A}(\bar{d})|\;[i,j]\subset I^{up}\}$.
\item[(ii)] $\mathcal{A}_+^{up}(\bar{d})=\{[i,j]\in
\mathcal{A}(\bar{d})|\;[i,j]\subset I_+^{up}\}$.
\item[(iii)] $\mathcal{A}^{down}(\bar{d})=\{[i,j]\in
\mathcal{A}(\bar{d})|\;[i,j]\subset I^{down}\}$.
\item[(iv)] $\mathcal{A}_+^{down}(\bar{d})=\{[i,j]\in
\mathcal{A}(\bar{d})|\;[i,j]\subset I_+^{down}\}$.
\end{itemize}
\end{definizione}
We restrict to $\Delta_{up}$ (one proceeds similarly for $\Delta_{down}$). Let $\bar{d}=\bar{d}_{up}+\bar{d}_{down}$ be a regular symplectic
dimension vector. $\Delta_{up}$
contains either a $\sigma_I$-fixed vertex, we shall call $n_{up}$, or a
$\sigma_I$-fixed edge, we shall call
$\xymatrix@-1pc{n_{up}\ar@{-}[r]&\sigma_I(n_{up})}$. Starting from
generic decomposition (\ref{dgvs}) of $\bar{d}_{up}$ we modify it
as follows.
\begin{itemize}
\item[(1)] We keep the summands $(e_{[i,j]}+\delta e_{[i,j]})^{\oplus q_{[i,j]}}$ corresponding to the
arc $[i,j]\subset I_+^{up}$.
\item[(2)] If $n_{up}$ is
such that there exists $a=a(n_{up})\in Q_1^{\sigma}$, then we keep
the summands $(e_{[i,\sigma_I(i)]})^{\oplus q_{[i,\sigma_I(i)]}}$
corresponding to the symmetric arcs $[i,\sigma_I(i)]$ of
$\Delta_{up}$.
\item[(3)] If $n_{up}$ is
such that there exists $x=x(n_{up})\in Q_0^{\sigma}$, we have the
symmetric dimension vectors
$$e_{[i_1,\sigma_I(i_1)]}
,\ldots,e_{[i_{2s},\sigma_I(i_{2s})]}$$ corresponding to the arcs
$[i_1,\sigma_I(i_1)],\ldots, [i_{2s},\sigma_I(i_{2s})]$ such that
$[i_1,\sigma_I(i_1)]\supseteq\cdots\supseteq
[i_{2s},\sigma_I(i_{2s})]$. Then we divide them into pairs
$$([i_{2k},\sigma_I(i_{2k})],[i_{2k-1},\sigma_I(i_{2k-1})]),$$ with
$1\leq k\leq s$. For each pair we consider
$[i_{2k},\sigma_I(i_{2k-1})]\cup [i_{2k-1},\sigma_I(i_{2k})]$ and
we substitute $e_{[i_{2k},\sigma_I(i_{2k})]}\oplus e_{
[i_{2k-1},\sigma_I(i_{2k-1})]}$ for
$$e_{[i_{2k},\sigma_I(i_{2k-1})]}+ e_{
[i_{2k-1},\sigma_I(i_{2k})]}.$$
\end{itemize}
So, by equation \ref{qi,j}, in the symplectic case we get (see lemma \ref{3218})
\begin{itemize}
\item[(i)] if $n_{up}$ is
such that there exists $a=a(n_{up})\in Q_1^{\sigma}$,
\begin{equation}\label{supqi,j1}
\bar{d}_{up}=\bigoplus_{[i,j]\in\mathcal{A}_+^{up}(\bar{d})
}(e_{[i,j]}+\delta e_{[i,j]})^{\oplus q_{[i,j]}} +
\bigoplus_{[i,\sigma_I(i)]\in\mathcal{A}^{up}(\bar{d})
}(e_{[i,\sigma_I(i)]})^{\oplus q_{[i,\sigma_I(i)]}};
\end{equation}
\item[(ii)] If $n_{up}$ is
such that there exists $x=x(n_{up})\in Q_0^{\sigma}$,
\begin{equation}\label{supqi,j2}
\bar{d}_{up}=\bigoplus_{[i,j]\in\mathcal{A}_+^{up}(\bar{d})
}(e_{[i,j]}+\delta e_{[i,j]})^{\oplus q_{[i,j]}}
+\bigoplus_{k=1}^{s}(e_{[i_{2k},\sigma_I(i_{2k-1})]}+ e_{
[i_{2k-1},\sigma_I(i_{2k})]}).
\end{equation}
\end{itemize}
Finally we have to modify like in (3) the dimension vector
$ph+e_{[i,\sigma_I(i)]}$ if $p$ is odd and $[i,\sigma_I(i)]$
passes through $n_{up}$ such that there exists $x=x(n_{up})\in
Q_0^{\sigma}$.
\begin{esempi}
Let $(Q,\sigma)$ be the symmetric quiver
$\widetilde{A}^{1,1}_{0,6}$. We recall that
$x_{\frac{l}{2}}=\sigma(x_{\frac{l}{2}})$. $\Delta$ has the form
(\ref{esempio}).\\
As in example \ref{esempioe}, let $p_1=4$, $p_2=3$, $p_3=0$ and
$p_4=2$. The $\sigma_I$-fixed vertex $4$ is such that
$e_4(x_{\frac{l}{2}})\neq 0$. The only symmetric arc passing
through 4 is $[4,4]$. Thus we substitute $(e_4)^{\oplus 2}$ for
$e_4+e_4$. So, in the symplectic case we get
$$
\sum_{i\in I_+}p_i(e_i+\delta e_i)+\sum_{i\in I_{\delta}}p_ie_i=
((e_2+\delta e_2)+e_1)^{\oplus 3}\oplus e_1\oplus 2
e_4.
$$
\end{esempi}
Similarly we proceed with the decomposition of
$\bar{d}'$ and $\bar{d}''$.\\
Let $\bar{d}=\bar{d}_{up}+\bar{d}_{down}$ be a regular orthogonal
dimension vector. Starting from generic
decomposition (\ref{dgvs}) of $\bar{d}_{up}$ we modify it as
follows.
\begin{itemize}
\item[(1)] We keep the summands $(e_{[i,j]}+\delta e_{[i,j]})^{\oplus q_{[i,j]}}$ corresponding to the
arc $[i,j]\subset I_+^{up}$.
\item[(2)] If $n_{up}$ is
such that there exists $a=a(n_{up})\in Q_1^{\sigma}$ such that
$\bar{d}_{ta}$ is odd or $n_{up}$ is such that there exist
$x=x(n_{up})\in Q_0^{\sigma}$, then we keep the summands
$(e_{[i,\sigma_I(i)]})^{\oplus q_{[i,\sigma_I(i)]}}$ corresponding
to the symmetric arcs $[i,\sigma_I(i)]$ of $\Delta_{up}$.
\item[(3)] If $n_{up}$ is
such that there exists $a=a(n_{up})\in Q_1^{\sigma}$ such that
$\bar{d}_{ta}$ is even, we have the symmetric dimension vectors
$$e_{[i_1,\sigma_I(i_1)]}
,\ldots,e_{[i_{2s},\sigma_I(i_{2s})]}$$ corresponding to the arcs
$[i_1,\sigma_I(i_1)],\ldots, [i_{2s},\sigma_I(i_{2s})]$ such that
$[i_1,\sigma_I(i_1)]\supseteq\cdots\supseteq
[i_{2s},\sigma_I(i_{2s})]$. Then we divide them into pairs
$$([i_{2k},\sigma_I(i_{2k})],[i_{2k-1},\sigma_I(i_{2k-1})]),$$ with
$1\leq k\leq s$. For each pair we consider
$[i_{2k},\sigma_I(i_{2k-1})]\cup [i_{2k-1},\sigma_I(i_{2k})]$ and
we substitute $e_{[i_{2k},\sigma_I(i_{2k})]}\oplus e_{
[i_{2k-1},\sigma_I(i_{2k-1})]}$ for
$$e_{[i_{2k},\sigma_I(i_{2k-1})]}+ e_{
[i_{2k-1},\sigma_I(i_{2k})]}.$$
\end{itemize}
So, by equation \ref{qi,j}, in the orthogonal case we get
\begin{itemize}
\item[(i)] if $n_{up}$ is
such that there exists $a=a(n_{up})\in Q_1^{\sigma}$ such that
$\bar{d}_{ta}$ is odd or $n_{up}$ is such that there exist
$x=x(n_{up})\in Q_0^{\sigma}$,
\begin{equation}\label{supqi,j1}
\bar{d}_{up}=\bigoplus_{[i,j]\in\mathcal{A}_+^{up}(d')
}(e_{[i,j]}+\delta e_{[i,j]})^{\oplus q_{[i,j]}} +
\bigoplus_{[i,\sigma_I(i)]\in\mathcal{A}^{up}(d')
}(e_{[i,\sigma_I(i)]})^{\oplus q_{[i,\sigma_I(i)]}};
\end{equation}
\item[(ii)] if $n_{up}$ is
such that there exists $a=a(n_{up})\in Q_1^{\sigma}$ such that
$\bar{d}_{ta}$ is even,
\begin{equation}\label{supqi,j2}
\bar{d}_{up}=\bigoplus_{[i,j]\in\mathcal{A}_+^{up}(d')
}(e_{[i,j]}+\delta e_{[i,j]})^{\oplus q_{[i,j]}}
+\bigoplus_{k=1}^{s}(e_{[i_{2k},\sigma_I(i_{2k-1})]}+ e_{
[i_{2k-1},\sigma_I(i_{2k})]}).
\end{equation}
\end{itemize}
Finally we have to modify like in (3) the dimension vector
$ph+e_{[i,\sigma_I(i)]}$ if $p$ is odd and $[i,\sigma_I(i)]$
passes through $n_{up}$ such that there exists $a=a(n_{up})\in
Q_1^{\sigma}$ such that $\bar{d}_{ta}$ is even.
\begin{esempi}
Let $(Q,\sigma)$ be the symmetric quiver
$\widetilde{A}^{1,1}_{0,6}$. We recall that $b=\sigma(b)$.
$\Delta$ has the form
(\ref{esempio}).\\
As in example \ref{esempioe}, let $p_1=4$, $p_2=3$, $p_3=0$ and
$p_4=2$. The $\sigma_I$-fixed vertex 1 is such that $e_1(tb)\neq
0$ and $\bar{d}_{tb}$ is 2. The only symmetric arcs passing through
1 is $[2,\sigma_I(2)]$. Thus we substitute $((e_2+\delta
e_2))+e_1)^{\oplus 3}\oplus e_1$ for $2((e_2+\delta e_2))+e_1)\oplus ((e_2+\delta e_2))+2e_1)$. So, in the
orthogonal case we get
$$
\sum_{i\in I_+}p_i(e_i+\delta e_i)+\sum_{i\in I_{\delta}}p_ie_i=
 2((e_2+\delta e_2)+
e_1)\oplus ((e_2+\delta e_2))+2e_1)\oplus (e_4)^{\oplus 2}.
$$
\end{esempi}
Similarly we proceed with the decomposition of $\bar{d}'$ and
$\bar{d}''$.\\
In general we have
\begin{proposizione}\label{dg}
Let $(Q,\sigma)$ be a symmetric quiver of tame type.
\begin{itemize}
\item[(1)] If $d$ is a regular symplectic dimension vector with decomposition
(\ref{fdc}). Then
\begin{equation}\label{fdg}
d=\bigoplus_{i=1}^p
h\oplus\bar{d}_{up}\oplus\bar{d}_{down}\oplus\bar{d}'_{up}\oplus\bar{d}'_{down}\oplus\bar{d}''_{up}\oplus\bar{d}''_{down}
\end{equation}
is the symplectic generic decomposition of $d$.
\item[(2)] If $d$ is a regular orthogonal dimension vector with decomposition
(\ref{fdc}). Then
\begin{equation}\label{fdgo}
d=\bigoplus_{i=1}^p
h\oplus\bar{d}_{up}\oplus\bar{d}_{down}\oplus\bar{d}'_{up}\oplus\bar{d}'_{down}\oplus\bar{d}''_{up}\oplus\bar{d}''_{down}
\end{equation}
is the orthogonal generic decomposition of $d$.
\end{itemize}
\end{proposizione}
For the proof, we need two propositions. We state these
propositions only for regular indecomposable symplectic
(respectively orthogonal) representations related to polygon
$\Delta$, because for those related to polygon $\Delta'$ and to
polygon $\Delta''$ the statement and the proof are similar.
\begin{proposizione}\label{extreg1}
Let $(Q,\sigma)$ be a symmetric quiver of tame tape. Let $V_1\neq
V_2$ be two regular indecomposable symplectic (respectively
orthogonal) representations of $(Q,\sigma)$ with symmetric
dimension vector corresponding respectively to the arc $[i,j]$ and
the arc $[h,k]$ of $\Delta$ ($\Delta'$ or $\Delta''$). Moreover we
suppose that $[i,j]$ and $[h,k]$ don't satisfy the following
properties
\begin{itemize}
\item[(i)] $[i,j]\cap [h,k]\neq \emptyset$ and $[i,j]$ doesn't
contain $[h,k]$;
\item[(ii)] $[i,j]\cap [h,k]\neq \emptyset$ and $[h,k]$ doesn't
contain $[i,j]$;
\item[(iii)] $[i,j]$ and $[h,k]$ are linked by one edge of
$\Delta$ (respectively $\Delta'$ or $\Delta''$).
\end{itemize}
Then $Ext^1_Q(V_1,V_2)=0$.
\end{proposizione}
\begin{proof}
We obtain the statement by
$$
\langle e_i,e_j\rangle=\left\{\begin{array}{ll}1&\textrm{if}\;i=j\\
-1&\textrm{if}\;i=j-1\\
0&\textrm{otherwise}.\end{array}\right.
$$ and applying Lemma 6.3 in chapter IX of \cite{ass} to the symplectic (respectively orthogonal) representations $V_1$ and $V_2$ corresponding to the arcs which don't satisfy (i), (ii) and (iii) (for more details, see Proposition 3.2.32 \cite{a1}).
\end{proof}
\begin{proposizione}\label{autoext}
Let $(Q,\sigma)$ be a symmetric quiver of tame tape. Let $V$ be a
regular indecomposable symplectic (respectively orthogonal)
representation of $(Q,\sigma)$ such that $\underline{dim}(V)=h$ or
$\bar{d}$. Moreover we suppose $V\neq E_{i,j}\oplus
E_{\sigma_I(j),\sigma_I(i)}$ with $i,j\in I_+$ such that
$e_i(ta)\neq 0$ or $e_j(ta)\neq 0$ for $a\in Q_1^{\sigma}$. Then,
for every non-trivial short exact sequence
$$
0\rightarrow V\rightarrow W\rightarrow V\rightarrow 0,
$$
$W$ is not symplectic (respectively it is not orthogonal).
\end{proposizione}
\begin{proof}
The statement follows from an analysis case by case of symmetric quivers of tame type describing the non-trivial autoextension of each regular indecomposable symplectic (respectively orthogonal) $V$ of dimension $h$ or $\bar{d}$ such that $V\neq E_{i,j}\oplus
E_{\sigma_I(j),\sigma_I(i)}$ with $i,j\in I_+$ such that
$e_i(ta)\neq 0$ or $e_j(ta)\neq 0$ for $a\in Q_1^{\sigma}$ (for details of the type $\widetilde{A}^{2,0,1}_{k,l}$, see Proposition 3.2.33 in \cite{a1}).
\end{proof}
\begin{proof}[Proof of Proposition \ref{dg}]
\textit{(1)} Let $d$ be a symplectic
regular dimension vector with decomposition (\ref{fdg}). First we
note that the symmetric dimension vectors appearing in
decomposition (\ref{fdg}) are not dimension vectors of the regular
indecomposable symplectic representations which are exceptions of
proposition \ref{extreg1} and \ref{autoext}. Let $\mathcal{O}(d)$
be the open orbit of the regular symplectic representations of
dimension $d$. By [Bo1] and [Z], we obtain each
representation $V$ in $\mathcal{O}(d)$ as follows.\\
There are representations $M_i$, $U_i$, $V_i$ and short exact
sequences
$$
0\rightarrow U_i\rightarrow M_i\rightarrow V_i\rightarrow 0
$$
such that $M_{i+1}=U_i\oplus V_i$ and $V=U_{n+1}\oplus V_{n+1}$,
with $1\leq i\leq n$ for some $n\in\mathbb{N}$.\\
By Propositions \ref{extreg1} and \ref{autoext}, we have
\begin{itemize}
\item[(i)] If $U_i\neq V_i$, then $Ext^1_Q(V_i,U_i)=0$.
\item[(ii)] If $U_i=V_i$, then either $Ext^1_Q(U_i,U_i)=0$ or
no one non-trivial auto-extension of $U_i$ is symplectic. So, if
$Ext^1_Q(U_i,U_i)\neq 0$ then $U_i$ doesn't appear in
decomposition of a symplectic representation.
\end{itemize}
Hence $V$ decomposes in regular indecomposable symplectic
representations of dimension $\beta_i$, where $\beta_i$ are
regular symmetric dimension
vectors appearing in decomposition (\ref{fdg}) of $d$.\\
\textit{(2)} One proves similarly to \textit{(1)}. 
\end{proof}

\subsubsection{Proof of Theorem \ref{stq}}
By
Proposition 2.4 in \cite{ar},  Proposition \ref{oA} and Lemma \ref{kac}, we can reduce the
proof to the orientation of $\widetilde{A}$ as in proposition
\ref{oA} and to the equiorientation for $\widetilde{D}$.\\
Let $d$ be a regular symmetric vector with a decomposition
(\ref{fdg}) or (\ref{fdgo}). We note that if $d=d_1\oplus d_2$ with
$d_1$ and $d_2$ summands of this generic decomposition, we have
canonical embeddings
\begin{equation}\label{Spe}
SpSI(Q,d)\stackrel{\Phi_d}{\rightarrow}\bigoplus_{\chi\in
char(Sp(Q,d))}SpSI(Q,d_1)_{\chi|_{d_1}}\otimes
SpSI(Q,d_2)_{\chi|_{d_2}}
\end{equation}
and
\begin{equation}\label{Oe}
OSI(Q,d)\stackrel{\Psi_d}{\rightarrow}\bigoplus_{\chi\in
char(O(Q,d))}OSI(Q,d_1)_{\chi|_{d_1}}\otimes
OSI(Q,d_2)_{\chi|_{d_2}},
\end{equation}
induced by the restriction homomorphism. We prove theorem
\ref{stq} by induction on the number of the summands
$e_{[i,j]}+\delta e_{[i,j]}$, $e_{[i,\sigma_I(i)]}$,
$e_{[i_{2k},\sigma_I(i_{2k-1})]}+e_{[i_{2k-1},\sigma_I(i_{2k})]}$
and respective summands corresponding to the admissible arcs in
$\mathcal{A}'(d)$ and in $\mathcal{A}''(d)$. If this number is 0,
then $d=ph$ and it was already proved. We suppose that the generic
decomposition of $d$ contains one of those summands and, without
loss of generality, we can assume that this summand is one of
those corresponding to the arcs in $\mathcal{A}(d)$. In particular
we suppose that this summand is $e_{[s,\sigma_I(s)]}$ (one
proceeds similarly for the other types), with $s\in I_+\sqcup
I_{\delta}$, and we can assume $ind[s,\sigma_I(s)]=r=max\{p_k\}$.
We call $d_2=e_{[s,\sigma_I(s)]}$ and so
$d_1=d-e_{[s,\sigma_I(s)]}$. Now we compare the generators of the
algebras $SpSI(Q,d)$ and $SpSI(Q,d_1)$ (respectively  $OSI(Q,d)$
and $OSI(Q,d_1)$). By induction the generators of $SpSI(Q,d_1)$
(respectively of $OSI(Q,d_1)$) are described by theorem \ref{stq}.
Since $\Delta'(d)=\Delta'(d_1)$ and $\Delta''(d)=\Delta''(d_1)$,
the generators $c_0,\ldots,c_t$ (with $t=\frac{p}{2}$,
$\frac{p-1}{2}$ or $p$), those corresponding to the arcs from
$\mathcal{A}'(d)$ and those corresponding to the arcs from
$\mathcal{A}''(d)$ occur. So it's enough to study the behavior of
the semi-invariants corresponding to the arcs from
$\mathcal{A}(d)$. We describe the link between the admissible arcs
of the polygons $\Delta(d)$ and $\Delta(d_1)$. We have
$$
d_1=ph+\sum_{i\in I_+\setminus(I_+\cap
[s,\sigma_I(s)])}p_i(e_i+\delta e_i)+\sum_{i\in
I_{\delta}\setminus(I_{\delta}\cap [s,\sigma_I(s)])}p_i e_i+
$$
$$
\sum_{i\in I_+\cap [s,\sigma_I(s)]}p_i(e_i+\delta e_i)+\sum_{i\in
I_{\delta}\cap [s,\sigma_I(s)]}p_i e_i+
$$
$$
\sum_{i\in I'_+}p'_i(e'_i+\delta e'_i)+\sum_{i\in I'_{\delta}}p'_i
e'_i+\sum_{i\in I''_+}p''_i(e''_i+\delta e''_i).
$$
We have two cases
\begin{itemize}
\item[(1)] $p_{s+1}=p_{\sigma_I(s)-1}<r-1$ with $s+1\in I_+$,
\item[(2)] $p_{s+1}=p_{\sigma_I(s)-1}=r-1$ with $s+1\in I_+$.
\end{itemize}
in the case (1) the only difference between the structure of
$\mathcal{A}(d)$ and $\mathcal{A}(d_1)$ is that the admissible
arcs $[s,s-1],[s-1,s-2],\ldots,[\sigma_I(s)+1,\sigma(s)]$ are of
index $r$ in $\mathcal{A}(d)$ and of index $r-1$ in
$\mathcal{A}(d_1)$. In the case (2) we have the admissible arc
$[s+1,\sigma_I(s)-1]$ of index $r-1$. The admissible arcs
$[s,s-1],[s-1,s-2],\ldots,[\sigma_I(s)+1,\sigma_I(s)]$ are of
index $s$ in $\mathcal{A}(d)$ and the admissible arcs
$[s+1,s],[s,s-1],\ldots,[\sigma_I(s)+1,\sigma_I(s)],[\sigma_I(s),\sigma_I(s)-1]$
are of index $r-1$ in $\mathcal{A}(d_1)$.\\
Now we prove that the embeddings $\Phi_d$ and $\Psi_d$ are
isomorphisms and this will be done in two steps. We describe the strategy for the symplectic case, similarly one proceeds for the orthogonal case. The first step is
to show case by case that the semi-invariants corresponding to the
admissible arcs $[i,j]$ are non zero $c^V$ for some $V\in Rep(Q)$
and, if $V$ satisfies property  \textit{(Op)}, they
are non zero $pf^V$. For the first step, we will deal only with the case $\widetilde{A}^{1,1}_{k,l}$, the other ones are similar. The second step is to give an explicit
description of the generators of the algebras on the right hand
side of $\Phi_d$ and $\Psi_d$. This is based on the knowledge,
given by inductive hypothesis, of the algebra $SpSI(Q,d_1)$. Moreover, by Proposition \ref{oa}, we know that $SpSI(Q,d_2)$ is a polynomial ring since $Sp(Q,d_2)$ has an open orbit in $SpRep(Q,d_2)$ (see Lemma \ref{sk}). So we can describe explicitly the generators of the algebra $SpSI(Q,d_2)$ and we note that they are determinants or pfaffians. At that point
we know the generators of the algebras on the right hand side of
$\Phi_d$ and $\Psi_d$. Now, using the fact that these are
determinants or pfaffians, we prove that they actually are in
$SpSI(Q,d)$ and that the
embeddings $\Phi_d$ and $\Psi_d$ are isomorphisms.\\
We will consider case by case the semi-invariants corresponding to
each admissible arc $[i,j]$ for $\widetilde{A}^{1,1}_{k,l}$ (for the other symmetric quivers of tame type see section 3.2.1 in \cite{a1}). To simplify the notation we shall
call $a$ both the arrow $a\in Q_1$ and the linear map $V(a)$
defined on $a$, where $V$ is a representation of $Q$.\\
We have at
most two $C^+$-orbits $\Delta$ and $\Delta'$ of the dimension
vectors of nonhomogeneous simple regular representation. We assume
$n\geq 2$ and we consider the $C^+$-orbit $$\{e_1=\delta
e_1,e_2,\ldots,e_{\frac{l}{2}},e_{\frac{l}{2}+1}=\delta
e_{\frac{l}{2}+1} ,\delta e_{\frac{l}{2}},\ldots,\delta
e_2\}.$$ Let $[i,j]\in\mathcal{A}(d)$. If we consider the arc
$[1,1]$ of index 0, i.e. $p_1=0,p_2\ne
0,\ldots,p_{\frac{l}{2}+1}\ne 0$, we have the minimal
projective resolution of $V_{(0,1)}$
$$
0\longrightarrow
P_{\sigma(a_0)}\stackrel{d^{V_{(0,1)}}_{min}}{\longrightarrow}
P_{a_0}\longrightarrow V_{(0,1)}\longrightarrow 0
$$
where
$d^{V_{(0,1)}}_{min}=\sigma(v_1)\cdots\sigma(v_{\frac{l}{2}})v_{\frac{l}{2}}\cdots
v_1$ and so
$$
c^{V_{(0,1)}}=det(Hom_Q(d^{V_{(0,1)}}_{min},\cdot))=det(\sigma(v_1)\cdots\sigma(v_{\frac{l}{2}})v_{\frac{l}{2}}\cdots
v_1)
$$
in the orthogonal case and
$pf^{V_{(0,1)}}=pf(\sigma(v_1)\cdots\sigma(v_{\frac{l}{2}})v_{\frac{l}{2}}\cdots
v_1)$ in the symplectic case, since by definition of symplectic
representation\\
$\sigma(v_1)\cdots\sigma(v_{\frac{l}{2}})v_{\frac{l}{2}}\cdots
v_1$ is skew-symmetric . If we consider the arc
$[2,\sigma_I(2)]=[2,0]$ of index 0, i.e.
$p_{\sigma_I(2)}=0=p_2,p_1\ne 0$, then we have the minimal
projective resolution of $V_{(1,0)}$
$$
0\longrightarrow
P_{\sigma(a_0)}\stackrel{d^{V_{(1,0)}}_{min}}{\longrightarrow}
P_{a_0}\longrightarrow V_{(1,0)}\longrightarrow 0
$$
where
$d^{V_{(1,0)}}_{min}=\sigma(u_1)\cdots\sigma(u_{\frac{k}{2}})bu_{\frac{k}{2}}\cdots
u_1$ and so
$$
c^{V_{(1,0)}}=det(Hom_Q(d^{V_{(1,0)}}_{min},\cdot))=
det(\sigma(u_1)\cdots\sigma(u_{\frac{k}{2}})bu_{\frac{k}{2}}\cdots
u_1)
$$
in the symplectic case and
$pf^{V_{(1,0)}}=pf(\sigma(u_1)\cdots\sigma(u_{\frac{k}{2}})bu_{\frac{k}{2}}\cdots
u_1)$ in the orthogonal case, since $b$ is skew-symmetric and
$\sigma(u_i)=-(u_i)^t$. We note that for $l=2$ we have only the
admissible arcs $[1,1]$ an $[2,\sigma_I(2)]$. We assume now that
$l\geq 4$ ($l$ is even) and $[i,j]$ is not an admissible arc
considered above. If $1\leq j\leq i-1\leq l$, then we identify $[i,j]$
with the path $v_{i-1}\cdots v_{j}$ in $Q$ and we have the minimal
projective resolution of $E_{i-1,j}$
$$
0\longrightarrow
P_{x_{i-1}}\stackrel{d^{E_{i-1,j}}_{min}}{\longrightarrow}
P_{x_{j-1}}\longrightarrow E_{i-1,j}\longrightarrow 0
$$
where $d^{E_{i-1,j}}_{min}=v_{i-1}\cdots v_{j}$ and so
$$
c^{E_{i-1,j}}=det(Hom_Q(d^{E_{i-1,j}}_{min},\cdot))=
det(v_{i-1}\cdots v_{j}).
$$
We note that
$$c^{\tau^-\nabla
E_{i-1,j}}=c^{E_{\sigma_I(j)-1,\sigma_I(i)}}=det(\sigma(v_{j})\cdots
\sigma(v_{i-1}))=det(v_{i-1}\cdots v_{j})=c^{E_{i-1,j}}.$$
Moreover, if $j=\sigma_I(i)$ then, only in the symplectic case, we
get $pf(\sigma(v_{i})\cdots v_{i})=pf^{E_{i-1,\sigma_I(i)}}$ since
$\sigma(v_{i})\cdots v_{i}$ is skew-symmetric. Now we consider the
arcs $[i,j]$ which have $e_1$ as internal vertex or as first vertex. For these arcs,
$1\leq i-1< j\leq l$ and we have the minimal projective resolution of
$E_{i-1,j}$
$$
0\longrightarrow P_{\sigma(a_0)}\oplus
P_{x_{i-1}}\stackrel{d^{E_{i-1,j}}_{min}}{\longrightarrow}
P_{a_0}\oplus P_{x_{j-1}}\longrightarrow E_{i-1,j}\longrightarrow
0
$$
where
$d^{E_{i-1,j}}_{min}=\left(\begin{array}{cc}\sigma(u_1)\cdots
b\cdots u_1
&\sigma(v_1)\cdots v_{j}\\
v_{i-1}\cdots v_1 & 0 \end{array}\right)$ and so
$$
c^{E_{i-1,j}}=det(Hom_Q(d^{E_{i-1,j}}_{min},\cdot))=det\left(\begin{array}{cc}\sigma(u_1)\cdots
b\cdots u_1
&v_{i-1}\cdots v_1\\
\sigma(v_1)\cdots v_{j} & 0 \end{array}\right).
$$
In particular we note that if $i=\sigma_I(j)$, in the orthogonal
case, we get
$$pf^{E_{\sigma_I(j)-1,j}}=pf\left(\begin{array}{cc}\sigma(u_1)\cdots
b\cdots u_1
&v_{i-1}\cdots v_1\\
\sigma(v_1)\cdots \sigma(v_{i-1}) & 0
\end{array}\right),$$ since $b$ is skew-symmetric, $\sigma(v_i)=-(v_i)^t$ and
 $\sigma(u_i)=-(u_i)^t$. Finally we note that
$V_{(0,1)}$, ${E_{i-1,\sigma_I(i)}}$ satisfy \textit{(Op)} and
$V_{(1,0)}$, ${E_{\sigma_I(j)-1,j}}$ satisfy property
\textit{(Spp)}. Similarly we define the semi-invariants for the
admissible arcs $[i,j]$ in $\mathcal{A}'(d)$, exchanging the upper
paths of $\widetilde{A}_{k,l}^{1,1}$ with the lower ones.\\
\\
We prove the second step of
proof of Theorem \ref{stq}. We note
that if $[i,j]$ is admissible then the semi-invariants associated
to $[i,j]$ define a nonzero element of $SpSI(Q,d)$ (respectively
of $OSI(Q,d)$).\\
For a symmetric dimension vector $d$ we denote
\begin{equation}\label{SpGamma}
Sp\Gamma(Q,d)=\{\chi\in\mathbb{Z}^{Q_0}\cup\frac{1}{2}\mathbb{Z}^{Q_0}|\,SpSI(Q,d)_{\chi}\neq
0\}
\end{equation}
and
\begin{equation}\label{OGamma}
O\Gamma(Q,d)=\{\chi\in\mathbb{Z}^{Q_0}\cup\frac{1}{2}\mathbb{Z}^{Q_0}|\,OSI(Q,d)_{\chi}\neq
0\}
\end{equation}
the semigroup of weights of symplectic (respectively orthogonal)
semi-invariants. We note that (\ref{SpGamma}) and (\ref{OGamma})
involve also $\frac{1}{2}\mathbb{Z}^{Q_0}$ because in $SpSI(Q,d)$
and in $OSI(Q,d)$ also pfaffians can appear. To simplify the
notation, we shall call $\chi_{[i,j]}$, $\chi'_{[i,j]}$ and
$\chi''_{[i,j]}$ be respectively the weights of the
semi-invariants associated to admissible arcs $[i,j]$ respectively
from $\mathcal{A}(d)$, $\mathcal{A}'(d)$ and $\mathcal{A}''(d)$.
In the next the following proposition will be useful. We will
state it only for $\Delta$, because for $\Delta'$ and $\Delta''$
the statements are similar. Let $d$ be a regular symmetric
dimension vector with decomposition $d=ph+d'$ with
$p\geq 1$.
\begin{proposizione}\label{pe}
Let $(Q,\sigma)$ be a symmetric quiver of tame type. Let $d_2$ be
of type $e_{[s,\sigma_I(s)]}$, $e_{[s,t]}+\delta e_{[s,t]}$ or
$e_{[i_{2k},\sigma_I(i_{2k-1})]}+e_{[i_{2k-1},\sigma_I(i_{2k})]}$.\\
(i) If $d_2=e_{[s,\sigma_I(s)]}$, then
\begin{itemize}
\item[(a)] For every arc $[i,j]$ of $\Delta'$ and $\Delta''$ we
have $\chi'_{[i,j]}|_{supp(d_2)},\chi''_{[i,j]}|_{supp(d_2)}\in
Sp\Gamma(Q,d_2)$ (respectively in $O\Gamma(Q,d_2)$).
\item[(b)] For every arc $[i,j]$ of $\Delta$ that doesn't intersect
 $[s,\sigma_I(s)]$ or contains $[s+1,\sigma_I(s)-1]$ we
have $\chi_{[i,j]}|_{supp(d_2)}\in Sp\Gamma(Q,d_2)$ (respectively
in $O\Gamma(Q,d_2)$).
\item[(c)] Let $\rho_1,\ldots,\rho_r$ be the weights of generators of the polynomial
algebra $SpSI(Q,d_2)$ (respectively $OSI(Q,d_2)$). Then $r\geq
n'-s$, where $n'\in I_+\sqcup I_{\delta}$ is either a
$\sigma_I$-fixed vertex or the extremal vertex of a
$\sigma_I$-fixed edge, and $\rho_1,\ldots,\rho_r$ can be reordered
such that
$\rho_1=\chi_{[s,s-1]},\ldots,\rho_{n'-s}=\chi_{[n'+1,n']}$ and
for every $m>n'-s$ we have $\langle\rho_m,e_n\rangle=0$ for
$n=s,\ldots,n'$.
\end{itemize}
(ii) Let $d_2=e_{[s,t]}+\delta e_{[s,t]}$, then
\begin{itemize}
\item[(a)] For every arc $[i,j]$ of $\Delta'$ and $\Delta''$ we
have $\chi'_{[i,j]}|_{supp(d_2)},\chi''_{[i,j]}|_{supp(d_2)}\in
Sp\Gamma(Q,d_2)$ (respectively in $O\Gamma(Q,d_2)$).
\item[(b)] For every symmetric arc $[i,j]$ of $\Delta$ that doesn't intersect
 $[s,t]\cup[\sigma_I(t),\sigma_I(s)]$ or contains $[s+1,\sigma_I(s+1)]$ or $[\sigma_I(t-1),t-1]$, we have
$\chi_{[i,j]}|_{supp(d_2)}\in Sp\Gamma(Q,d_2)$ (respectively in
$O\Gamma(Q,d_2)$).
\item[(c)] For every arc $[i,j]\subset I_+$ (respectively
$[i,j]\subset I_-$) that doesn't intersect $[s,t]$ (respectively
$[\sigma_I(t),\sigma_I(s)]$) or contains $[s+1,t-1]$ we have
$\chi_{[i,j]}|_{supp(d_2)}\in Sp\Gamma(Q,d_2)$ (respectively in
$O\Gamma(Q,d_2)$).
\item[(d)] Let $\rho_1,\ldots,\rho_r$ be the weights of generators of the polynomial
algebra $SpSI(Q,d_2)$ (respectively $OSI(Q,d_2)$). Then $r\geq
t-s$ and $\rho_1,\ldots,\rho_r$ can be reordered such that
$\rho_1=\chi_{[s,s-1]},\ldots,\rho_{t-s}=\chi_{[t+1,t]}$ and for
every $m>t-s$ we have $\langle\rho_m,e_n\rangle=0$ for
$n=s,\ldots,t$.
\end{itemize}
(iii) Let
$d_2=e_{[i_{2k},i_{\sigma_I(i_{2k-1})}]}+e_{[i_{2k-1},i_{\sigma_I(i_{2k})}]}$,
then
\begin{itemize}
\item[(a)] For every arc $[i,j]$ of $\Delta'$ and $\Delta''$ we
have $\chi'_{[i,j]}|_{supp(d_2)},\chi''_{[i,j]}|_{supp(d_2)}\in
Sp\Gamma(Q,d_2)$ (respectively in $O\Gamma(Q,d_2)$).
\item[(b)] For every arc $[i,j]$ of $\Delta$ that doesn't intersect
 $[i_{2k-1},\sigma_I(i_{2k-1})]$ or contains $[i_{2k-1}+1,\sigma_I(i_{2k-1})-1]$ we
have $\chi_{[i,j]}|_{supp(d_2)}\in Sp\Gamma(Q,d_2)$ (respectively
in $O\Gamma(Q,d_2)$).
\item[(c)] Let $\rho_1,\ldots,\rho_r$ be the weights of generators of the polynomial
algebra $SpSI(Q,d_2)$ (respectively $OSI(Q,d_2)$). Then $r\geq
n'-s$, where $n'\in I_+\sqcup I_{\delta}$ is either a
$\sigma_I$-fixed vertex or the extremal vertex of a
$\sigma_I$-fixed edge, and $\rho_1,\ldots,\rho_r$ can be reordered
such that
$\rho_1=\chi_{[s,s-1]},\ldots,\rho_{n'-s}=\chi_{[n'+1,n']}$ and
for every $m>n'-s$ we have $\langle\rho_m,e_n\rangle=0$ for
$n=s,\ldots,n'$.
\end{itemize}
\end{proposizione}
\begin{proof} It proceeds type by type analysis. We prove only the symplectic case for
$Q=\widetilde{A}^{1,1}_{k,l}$ and for $d_2=e_{[s,\sigma_I(s)]}$,
because the procedure to prove all other cases is similar. We
order the vertices of $\widetilde{A}^{1,1}_{k,l}$ such that the
only source is 1 (so the only sink is $\sigma(1)$), $hv_{i-1}=i$
for every $i\in\{2,\ldots,\frac{l}{2}+1\}$, $hu_i=\frac{l}{2}+i+1$
for every $i\in\{1,\ldots,\frac{k}{2}\}$ and then the respective
conjugates by $\sigma$ of these. We shall call
$w_{(t^1)_{i_1},\ldots,(t^f)_{i_f}}$, where
$t^1,\ldots,t^f\in\mathbb{Z}\cup\frac{1}{2}\mathbb{Z}$ and
$\{i_1,\ldots,i_f\}$ is an ordered subset of
$\{1,\ldots,\frac{l}{2}+\frac{k}{2}+1,\sigma(\frac{l}{2}+\frac{k}{2}+1)=\frac{l}{2}+\frac{k}{2}+2,\ldots,\sigma(1)=l+k+1\}$,
the vector such that
$$
w_{(t^1)_{i_1},\ldots,(t^f)_{i_f}}(y)=\left\{\begin{array}{ll}(t^j)_{i_j}
& y=i_j,\forall j=1,\ldots,f\\
0&\textrm{otherwise}.\end{array}\right.
$$
Moreover we can associate in bijective way the vertex
$i\in\{2,\ldots,\frac{l}{2}\}\subset(\widetilde{A}^{1,1}_{k,l})_0^+$
to $i\in I_+$, the vertex $\frac{l}{2}+i+1$ of
$\widetilde{A}^{1,1}_{k,l}$ to $i+1\in I'_+$ and
the vertex $\frac{l}{2}+1$ to $\frac{l}{2}+1\in I_{\delta}$.\\
\textit{(a)} We have
$$
\chi'_{[i,j]}=w_{(1)_{\frac{l}{2}+j},(-1)_{\frac{l}{2}+i}}\quad\textrm{for}\quad
1\leq j\leq i-1\leq k+1,
$$
if $[i,j]$ has not $e_1$ as internal vertex;
$$
\chi'_{[i,j]}=w_{(1)_1,(-1)_{\frac{l}{2}+i},(1)_{\frac{l}{2}+j},(-1)_{\sigma(1)}}\quad\textrm{for}\quad
1\leq i<j-1\leq k+1
$$
if $[i,j]$ has $e_1$ as internal vertex and in particular if
$j=\sigma_I(i)$ we have
$$
\chi'_{[i,j]}=w_{(\frac{1}{2})_1,(-\frac{1}{2})_{\frac{l}{2}+i},(\frac{1}{2})_{\sigma(\frac{l}{2}+i)},(-\frac{1}{2})_{\sigma(1)}}.
$$
Now if $\langle\chi'_{[i,j]},e_{[s,\sigma_I(s)]}\rangle\neq 0$
then $\chi'_{[i,j]}\not\in SpSI(Q,d_2)$, but we note that
$\langle\chi'_{[i,j]},e_{[s,\sigma_I(s)]}\rangle=0$ for every $i$
and $j$, so we have \textit{(a)}.\\
\textit{(b)} We have
$$
\chi_{[i,j]}=w_{(1)_{j},(-1)_{i}}\quad\textrm{for}\quad 1\leq
j\leq i-1\leq
l\quad\textrm{and}\quad\chi_{[\sigma(i),i]}=w_{(\frac{1}{2})_{i},(-\frac{1}{2})_{\sigma(i)}}
$$
if $[i,j]$ has not $e_1$ as internal vertex or as first vertex;
$$
\chi_{[i,j]}=w_{(1)_1,(-1)_{i},(1)_{j},(-1)_{\sigma(1)}}\quad\textrm{for}\quad
1\leq i-1<j\leq l.
$$
if $[i,j]$ has $e_1$ as internal vertex or as first vertex.\\
Now we note that
$\langle\chi_{[i,j]},e_{[s,\sigma_I(s)]}\rangle\neq 0$ if
$[i,j]\cap[s,t]\neq\emptyset$ and
$[i,j]\nsupseteq[s-1,\sigma(s+1)=\sigma(s)-1]$, so we have
\textit{(b)}.\\
\textit{(c)} First we note that we can choose symmetric arcs of
each length from a fixed vertex of $\Delta$, because the result of
Theorem \ref{stq} is invariant respect to the Coxeter
transformation $\tau^+$. The generators of $SpSI(Q,d_2)$ associated to
$\Delta(d_2)$ are $c^{E_i}=det(v_i)$ of weight
$\chi_{[i,i+1]}=w_{(1)_{i},(-1)_{i+1}}$ for every
$i\in\{1,\ldots,s-1\}$ and
$c^{E_{s-1,\sigma_I(s)}}=det\left(\begin{array}{cc}\sigma(u_1)\cdots
b\cdots u_1
&v_{s}\cdots v_1\\
 \sigma(v_1)\cdots \sigma(v_{s})& 0 \end{array}\right)$ of weight
$\chi_{[s,\sigma_I(s)]}=w_{(1)_{1},(-1)_{s},(1)_{\sigma(s)},(-1)_{\sigma(1)}}$.
So we call $\rho_i=\chi_{[i,i+1]}$  for every
$i\in\{1,\ldots,s-1\}$ and $\rho_{n'-s}=\chi_{[s,\sigma_I(s)]}$,
where in this case $n'=\frac{l}{2}+1$. The other generators
are associated to $\Delta'(d_2)$ and so, as done in the part
\textit{(a)} of this proposition, their weight $\rho_m$, for
$m\in\{n'-s+1,\ldots,r\}$, are such that
$\langle\rho_m,e_n\rangle=0$ for
$n\in\{s,\ldots,n'\}$. 
\end{proof}
We assume now that $d=d_1+d_2$ where $d_1=ph+d'_1$ with $p\geq 1$
and $d_2=e_{[s,\sigma_I(s)]}$, $e_{[s,t]}+\delta e_{s,t]}$ or
$e_{[i_{2k},i_{\sigma_I(i_{2k-1})}]}+e_{[i_{2k-1},i_{\sigma_I(i_{2k})}]}$.
So we take the corresponding arc in a chosen position (for which
we proved proposition \ref{pe}).
\begin{proposizione}\label{d1d2}
Let $d,d_1,d_2$ be as above. We suppose that the semigroup
$Sp\Gamma(Q,d_1)$ (respectively $O\Gamma(Q,d_1)$) is generated by
the weights $\chi_{[i,j]}$, $\chi'_{[i,j]}$, $\chi''_{[i,j]}$ for
admissible arcs $[i,j]$ of the labelled polygons $\Delta(d_1)$,
$\Delta'(d_1)$, $\Delta''(d_1)$. Then $Sp\Gamma(Q,d_1)\cap
Sp\Gamma(Q,d_2)$ (respectively $O\Gamma(Q,d_1)\cap
O\Gamma(Q,d_2)$) is generated by the weights $\chi_{[i,j]}$,
$\chi'_{[i,j]}$, $\chi''_{[i,j]}$ for admissible arcs $[i,j]$ of
the labelled polygons $\Delta(d)$, $\Delta'(d)$, $\Delta''(d)$.
\end{proposizione}
\begin{proof} 
We prove it only for the othogonal case and for $d_2=e_{[s,\sigma_I(s)]}$, because for the symplectic case e the other types of $d_2$, the proof is similar.\\
We are two cases.\\
(1) Assume $p_{s+1}=p_{\sigma_I(s)-1}<r-1$. The admissible arcs of
$\Delta(d_1)$, $\Delta'(d_1)$, $\Delta''(d_1)$ and $\Delta(d)$,
$\Delta'(d)$, $\Delta''(d)$ are the same. By Proposition \ref{pe}
$O\Gamma(Q,d_2)$ contains
$\chi_{[s,s-1]},\ldots,\chi_{[\sigma_I(s)+1,\sigma_I(s)]}$ and all
the other weights corresponding to the admissible arcs of
$\Delta(d)$, $\Delta'(d)$ and $\Delta''(d)$.\\
(2) Assume $p_{s+1}=p_{\sigma_I(s)-1}=r-1$. We prove that
$O\Gamma(Q,d_1)\cap O\Gamma(Q, d_2)$ is generated by
$\chi'_{[i,j]}$ for every admissible arc $[i,j]$ of
$\Delta'(d_1)=\Delta'(d)$, $\chi''_{[i,j]}$ for every admissible
arc $[i,j]$ of $\Delta''(d_1)=\Delta''(d)$ and $\chi_{[i,j]}$ for
every admissible arc $[i,j]$ of $\Delta(d_1)$ of index smaller
than $r-1$ or not intersecting $[s,\sigma_I(s)]$, i.e.
$\chi_{[s,s-1]},\ldots,\chi_{[\sigma_I(s)+1,\sigma_I(s)]}$ and
$\chi_{[s+1,\sigma_I(s)-1]}=\chi_{[s+1,s]}+\cdots+\chi_{[\sigma_I(s),\sigma_I(s)-1]}$.
Let
$$
\chi=\sum_{[i,j]\in\mathcal{A}(d_1)}n_{i,j}\chi_{[i,j]}+\sum_{[i,j]\in\mathcal{A}'(d_1)}n'_{i,j}\chi'_{[i,j]}
+\sum_{[i,j]\in\mathcal{A}''(d_1)}n''_{i,j}\chi''_{[i,j]},
$$
with $n_{i,j},n'_{i,j},n''_{i,j}\geq 0$, be an element of
$O\Gamma(Q,d_1)$. We assume that $\chi$ is also in
$O\Gamma(Q,d_2)$. By Proposition \ref{pe}, we note that all the
generators of $O\Gamma(Q,d_1)$ except of $\chi_{[s+1,s]}$ and
$\chi_{[\sigma_I(s),\sigma_I(s)-1]}$ are also in $O\Gamma(Q,d_2)$.
Hence, if $\chi$ contains neither $\chi_{[s+1,s]}$ nor
$\chi_{[\sigma_I(s),\sigma_I(s)-1]}$, then $\chi$ is a linear
combination of desired generators. So we have to prove that if
$\chi$ contains $\chi_{[s+1,s]}$ (resp.
$\chi_{[\sigma_I(s),\sigma_I(s)-1]}$) with positive coefficient,
then it contains
$\chi_{[s,s-1]},\ldots,\chi_{[\sigma_I(s),\sigma_I(s)-1]}$ (resp.
$\chi_{[s+1,s]},\ldots,\chi_{[\sigma_I(s)+1,\sigma_I(s)]}$). Thus
we can subtract $\chi_{[s+1,\sigma_I(s)-1]}$ from
$\chi$.\\
We assume that $\chi$ contains $\chi_{[s+1,s]}$ with positive
coefficient (the proof is similar for
$\chi_{[\sigma_I(s),\sigma_I(s)-1]}$). We note that
$\langle\chi_{[s+1,s]},e_s\rangle=1$ and, by Proposition
\ref{pe}, the other generators of $O\Gamma(Q,d_1)$, except
$\chi_{[s,s-1]}$, have zero product scalar with $e_s$. Moreover,
$\chi\in O\Gamma(Q,d_2)$ and so, by Proposition \ref{pe},
$\langle\chi,e_s\rangle\leq 0$. Hence $\chi$ contains
$\chi_{[s,s-1]}$ with positive coefficient. By Proposition
\ref{pe}, it follows that $\langle\chi,e_s+e_{s-1}\rangle\leq 0$.
But $\langle\chi_{[s+1,s]}+\chi_{[s,s-1]},e_s+e_{s-1}\rangle=1$
and $\chi_{[s-1,s-2]}$ is the only generator of $O\Gamma(Q,d_1)$
with negative scalar product with $e_s+e_{s-1}$. Continuing in
this way, we check that $\chi$ contains
$\chi_{[s+1,s]},\chi_{[s,s-1]},\ldots,\chi_{[\sigma_I(s)+1,\sigma_I(s)]},\chi_{[\sigma_I(s),\sigma_I(s)-1]}$
with positive coefficients. So we can subtract
$\chi_{[s+1.\sigma_I(s)-1]}$ from $\chi$ and continue. In this way
we complete the proof. 
\end{proof}
Now we can finish the proof of Theorem \ref{stq}. Since  Theorem
\ref{stq} is equivalent to Theorems \ref{tp1} and \ref{tp2} for
tame type and regular dimension vectors, then, in this way, we
finish also the proof of theorems \ref{tp1} and \ref{tp2}.\\
Again we consider the embeddings
\begin{equation}\label{Spe1}
SpSI(Q,d)\stackrel{\Phi_d}{\rightarrow}\bigoplus_{\chi\in
char(Sp(Q,d))}SpSI(Q,d_1)_{\chi|_{d_1}}\otimes
SpSI(Q,d_2)_{\chi|_{d_2}}
\end{equation}
and
\begin{equation}\label{Oe1}
OSI(Q,d)\stackrel{\Psi_d}{\rightarrow}\bigoplus_{\chi\in
char(O(Q,d))}OSI(Q,d_1)_{\chi|_{d_1}}\otimes
OSI(Q,d_2)_{\chi|_{d_2}}
\end{equation}
where $Q$, $d$, $d_1$ and $d_2$ are as above. The semigroup of
weights of the right hand side of $\Phi_d$ and $\Psi_d$ are
respectively $Sp\Gamma(Q,d_1)\cap Sp\Gamma(Q,d_2)$ and
$O\Gamma(Q,d_1)\cap O\Gamma(Q,d_2)$. These are generated by
$\chi_{[i,j]}$, $\chi'_{[i,j]}$, $\chi''_{[i,j]}$ for admissible
arcs $[i,j]$ of the labelled polygons $\Delta(d)$, $\Delta'(d)$,
$\Delta''(d)$, by proposition \ref{d1d2}. So the algebras on the
right hand side of $\Phi_d$ and $\Psi_d$ are generated by the
semi-invariants of weights $\chi_{[i,j]}$, $\chi'_{[i,j]}$,
$\chi''_{[i,j]}$ and by the semi-invariants of weights $\langle
h,\cdot\rangle$ (or $\frac{1}{2}\langle h,\cdot\rangle$).\\
Finally, we note that the embeddings $\Phi_d$ and $\Psi_d$ are
isomorphisms because they are also isomorphisms in the weight
$\langle h,\cdot\rangle$ (or $\frac{1}{2}\langle h,\cdot\rangle$)
and so we completed the proof of Theorem \ref{stq}. Moreover, in
that way, we also proved Proposition \ref{fij}, expliciting the
semi-invariants of type $c^V$ for every admissible arc $[i,j]$,
and the following theorem
\begin{teorema}\label{dtso} Let $(Q,\sigma)$ be a symmetric quiver of
tame type and let $d$ be a regular symmetric dimension vector. We consider decomposition (\ref{fdc}) with $p\geq 1$ and $d'\neq 0$. There exist
isomorphisms of algebras
\begin{equation}SpSI(Q,d)\stackrel{\Phi_d}{\rightarrow}\bigoplus_{\chi\in
char(Sp(Q,d))}SpSI(Q,ph)_{\chi}\otimes
SpSI(Q,d')_{\chi'}\end{equation} and
\begin{equation}OSI(Q,d)\stackrel{\Psi_d}{\rightarrow}\bigoplus_{\chi\in
char(O(Q,d))}OSI(Q,ph)_{\chi}\otimes
OSI(Q,d')_{\chi'},\end{equation}
where $\chi'=\chi|_{d'}$, i.e.
the restriction of the weight $\chi$ to the support of $d'$.
\end{teorema}

\end{document}